\documentclass[10pt]{article}
\usepackage{amsmath,amsfonts,amsthm,amssymb,amscd}
\numberwithin{equation}{section}

\newcommand{\be}{\begin{equation}}
\newcommand{\ee}{\end{equation}}
\newcommand{\bs}{\begin{split}}

\newcommand{\es}{\end{split}}
\newcommand{\ba}{\begin{align}}
\newcommand{\ea}{\end{align}}
\newcommand{\basl}[1]{\begin{align}\begin{split}\label{#1}}
\newcommand{\bas}{\begin{align}\begin{split}}

\newtheorem{theo}{Theorem}[section]
\newtheorem{prop}[theo]{Proposition}
\newtheorem{lemm}[theo]{Lemma}

\newtheorem{defi}[theo]{Definition}

\newcommand\fpr{\hfill$\Box$\null}

\newcommand\R{\mathbb{R}}
\newcommand\C{\mathbb{C}}

\title{Infinite dimensional semiclassical analysis and applications\\ to a model in NMR}

\author{ L. Amour, L. Jager, J. Nourrigat}

\date{Universit\'e de Reims, France}

\begin{document}

\maketitle

\begin{abstract}
\noindent
We are interested in this paper with the connection between the dynamics of a model related to Nuclear Magnetic Resonance (NMR) in  Quantum Field Theory (QFT)  with its classical counterpart known as the Maxwell-Bloch equations.
The model in QFT is a model of Quantum Electrodynamics (QED) considering fixed spins interacting with the quantized electromagnetic field in an external constant magnetic field. This model is close to the common spin-boson model. The classical model goes back to F. Bloch \cite{BL} in 1946. 
Our goal is not only to study the derivation of the Maxwell-Bloch equations  but to also  establish a semiclassical asymptotic expansion of  arbitrary high orders with control of the error terms of this standard nonlinear classical motion equations. This provides therefore  quantum 
 corrections of any order in powers of the semiclassical parameter of the Bloch equations. Besides, the asymptotic expansion for the photon number is also analyzed and a law describing the photon number time evolution is written down involving the radiation field polarization.
Since the quantum photon state Hilbert space (radiation field)  is infinite dimensional we are thus concerned in this article with the issue of semiclassical calculus in an infinite dimensional setting. In this regard, we are studying standard notions as Wick and anti-Wick quantizations,  heat operator, Beals characterization theorem and compositions of symbols in the infinite dimensional context which can have their own interest.
\end{abstract}

\parindent=0pt

\

{\it Keywords:} Semiclassical analysis, infinite dimensional analysis, composition of operators, Wick quantization, anti-Wick quantization, Wick symbol, Husimi function, Wiener spaces, Heat operator, symbolic calculus, QED, quantum electrodynamics, Maxwell-Bloch equations, Bloch equations, NMR, Nuclear Magnetic Resonance, photon emission, photon number.

\

{\it MSC 2010:} 35S05, 81V10, 47G30, 81Q20, 47L80, 26E15, 28C20.

\tableofcontents

\parindent = 0 cm

\parskip 10pt
\baselineskip 15 pt

\section{Statement of the results.}\label{s-1}

The Wick, Weyl and anti-Wick quantizations when  they depend on a
semiclassical parameter $h>0$ can be used, among other things,
to establish a link between quantum observable time evolutions
 of a quantum system with the classical system corresponding evolutions. For that purpose, one can determine a semiclassical expansion of
one of the symbols of the quantum time evolving observable. In a finite dimensional
setting, this is a consequence of particular cases of Egorov Theorem
(see \cite{Rob}, \cite{Zwo}, $\dots$).
This connection between both dynamics can be brought to light with the help of
coherent states, see \cite{CR} and \cite{Hepp},
or by constructing semiclassical measures, see  \cite{Ammari-Falconi} and
  \cite{A-N}.
 However, one notes that 
 symbolic calculus provides an asymptotic expansion
of an arbitrarily high order and not only the first term of the expansion.

The aim of this work is then to carry on the study of a symbolic calculus
and to apply it to the semiclassical limit of the evolution for a  quantum field model
in Nuclear Magnetic Resonance (NMR), see Section 4.11 in  \cite{Reu}, \cite{E} and \cite{J-H}.
In this perspective, we are here interested in the interaction between $N$  spin-$\frac{1}{2}$ particles   with the quantized electromagnetic field together with a constant external magnetic field. The spin particles are fixed 
at some arbitrary points
 $x_{\lambda }$ of $\R^3$.
This interaction model is also closely related to the spin-boson model (for example see,
\cite{H-S}, \cite{A-H}, $\dots$).

Let us first recall some usual models describing the same physical model.

The simplest model goes back to F. Bloch \cite{BL} in 1946. In this first model, spins are viewed as vectors ${\bf S}^{\lambda } (t)\in \R^3$, $\lambda=1,\dots,N$. Time evolutions of these vectors follow the Bloch equations (1946):
 $$\frac {d} {dt}  {\bf S}  ^{\lambda } (t) = 2 {\bf B} ^{ext } (x_{\lambda} , t) )
\times {\bf S}  ^{\lambda } (t)$$
where ${\bf B} ^{ext } (x, t)$ is a non quantized external magnetic  field. The spins affect the non quantized magnetic field according  to Maxwell equations through a current density  ${\bf j}(x , t)$.

One handicap of this model is the long time behavior of the solutions. Indeed,  the behavior of the solutions is not always physically consistent for large time $t$. To overcome that difficulty, one usually inserts  ad hoc additional terms in the Bloch equations \cite{BL} in the classical approach,
where the spins of the particles are still considered as vectors of  $\R^3$ with 
evolutions described by the Bloch equations (see (\ref{9-3})) and with 
 electric and magnetic fields still evolving according to the Maxwell equations.

A second model can also be considered. The spins are now quantized but the electromagnetic field is again not quantized. The configuration space of the system of $N$ spins is then the space ${\cal H}_{sp} =( \C^2  )^{\otimes N} $.
The fermion property for the spin-$\frac{1}{2}$ fixed particles is omitted here.
  In the space  ${\cal H}_{sp} $, we shall use the operators related to the
spins of the different particles. Let  $\sigma _j$ ($1 \leq j \leq 3$)
be the Pauli matrices:
\be\label{Pauli} \sigma_1 = \begin{pmatrix}  0 & 1 \\ 1 & 0   \end{pmatrix},
\qquad
\sigma_2 = \begin{pmatrix} 0 & -i \\ i & 0   \end{pmatrix},
\qquad
\sigma_3 = \begin{pmatrix}  1 & 0 \\ 0 & -1  \end{pmatrix}.\ee
For all $\lambda \leq N$ and all  $m\leq 3$, we denote by $\sigma_m^{[\lambda]}$
the operator in ${\cal H} _{sp}$ defined by:
\be\label{spin-ini}\sigma_m^{[\lambda]} = I \otimes \cdots  \otimes I \otimes \sigma_m\otimes  I  \otimes \cdots \otimes I,
\ee
where  $\sigma_m$ is located at the  $\lambda ^{th}$ position.
 Namely, the Hamiltonian of the spin system in an external magnetic field ${\bf B}^{ext}(x)$ is:
$$ H_{mag} = h \sum _{\lambda =1}^N  \sum _{m=1}^3  B_m^{ext} (x_{\lambda}) \otimes
\sigma_m^{[\lambda]}. $$
Then, the time evolution $\psi (t)$ of an initial state $\psi$ belonging to ${\cal H}_{sp}$ is defined as:
$$ i h \frac {\partial \psi(t)} {\partial t} = H_{mag}\psi (t),\qquad \psi(0)=\psi. $$
This model is commonly used in quantum computing (see \cite{VdS-C}). Elementary operations in quantum computing are called "quantum gates" and are unitary operators in ${\cal H}_{sp}$. An artificial interaction between the spins is needed.
Spins are thus interacting with a non quantized external magnetic field without modifying the field.

In the third model, the spins and the electromagnetic field are quantized. This is the model  under consideration in this work. The Hilbert space and Hamiltonian operator 
are depending on the semiclassical parameter  $h>0$ and are described below in this subsection.  Let us also mention that, this model  seems moreover necessary in order to study the part played  by possible electrons in the spin-nucleus interaction. See \cite{Ro-Au} for more details concerning this last point. The motion of the electrons is not taken into account in our work.

Our purpose here is to prove that the first model is a semiclassical limit of the third model (see (\ref{MaxSym-5})(\ref{MaxSym-6}) and also (\ref{9-3})(\ref{9-4})). This can therefore be viewed as the derivation of the Bloch equations. Note that a full asymptotic expansion with control of the error terms is in addition obtained. As already said, these semiclassical approximations concern observable time evolutions. We point out that \cite{J-H} (see also \cite{E}) adresses a similar issue.

We now give the third model written in the context of QED.

The one photon configuration Hilbert  space $H$  is the set of mappings  $f \in L^2 (\R^3, \R^3)$ satisfying $k\cdot f(k) = 0$ almost everywhere in $k\in\R^3$ (see \cite{L-L}) where
$|f|^2=\int_{\R^3}|f(k)|^2dk$. The Hilbert space ${\cal H}_{ph}$ of photon quantum states is
 the symmetrized Fock space
${\cal F}_s (H_{\bf C})$ over  the complexified of $H$. 
We follow \cite{RSII} for Fock spaces considerations and notations, in particular, for the usual operators in these spaces:
the Segal field $\Phi_S(V)$ associated with an element $V$ in
$H^2$, the $\Gamma (T)$ and ${\rm d} \Gamma (T)$ operators associated with some operators
$T$ acting in $H^2$. Note that, throughout this paper, the space $H^2$ is sometimes identified with the complexified
$H_{\bf C}$ but this identification is not everywhere systematically effectuated in order to avoid possible confusions. 
With the aim of underlining the role of the semiclassical parameter
 $h>0$, we also set for all $V$ belonging to $H^2$:
\be\label{Phi-S-h} \Phi_{S h} (V)=  h^\frac{1}{2} \Phi_S(V).\ee

The Hilbert space describing the states of $N$ fixed particles with spin-$\frac{1}{2}$,
without interaction, at a given time, is  ${\cal H}_{sp}  =( \C^2  )^{\otimes N}$.
Again, the fermionic property of the spin-$\frac{1}{2}$ fixed particles is not taken into account here. The Hilbert space describing the states of the whole system, consisting of
the quantum field and of the $N$ particles, in the presence of a constant
magnetic field $\beta$, is the completed tensor product
${\cal H}_{ph} \otimes {\cal H}_{sp} $.

Let  $M_{\omega}$ be the operator with domain   $D(M_{\omega}) \subset H$ such
that $M_{\omega} q (k) = |k| q(k)$ almost everywhere in $k\in\R^3$.
We denote in the same way the analogous operators defined  on $H^2$ or
on the complexified $H_{\bf C}$. In the Fock space framework, the photon free energy Hamiltonian operator
$H_{ph}$  is usually defined as  $h {\rm d} \Gamma (M_{\omega})$.

The photon number operator denoted by $N$ is $N = {\rm d} \Gamma (I)$.

The three components of the magnetic field at each point $x$ in  $\R^3$ are defined using the elements 
 $ B_{jx}$ belonging to  $H^2$ and written as
follows,  when one identifies $H^2$ with the complexified $H_{\bf C}$:
\be\label{7.3} B_{jx}(k) = {i\chi(|k|)|k|^{1\over 2} \over (2\pi)^{3\over 2}}
e^{-i(   k\cdot x   )} {k\times e_j \over |k|},\quad k\in\R^3\backslash\{0\}\ee
where the function $\chi $ (ultraviolet cut-off) belongs to ${\cal S} (\R)$.
Also set,
\be\label{7.4} E_{mx}   =  J B_{mx},\quad  \leq m \leq 3 \ee
where  $J: H^2 \rightarrow H^2$ stands for the helicity operator defined by, 
\be\label{7.5}JX  (k) =  {k\times X(k) \over |k| },
\quad k \in \R^3 \setminus \{ 0 \},\ee
for each $X$ in $H^2$
or $H_{\bf C}$.
One then defines the electric and magnetic fields components operators at each point 
 $x$ of $\R^3$ by:
$$B_m(x) = \Phi_{S h} (B_{mx})=  {h}^\frac{1}{2} \Phi_S(B_{mx}),  \quad E_m(x) =
\Phi_{S h}(E_{mx})=  {h}^\frac{1}{2} \Phi_S(E_{mx}),$$
for $m=1,2,3$.

{\it The Hamiltonian.} The Hamiltonian of the system studied here is the Pauli Fierz Hamiltonian where the terms concerning the spin particles motion are deleted. This Hamiltonian is often used for modeling NMR in quantum field theory (see  \cite{Ro-Au}\cite{E}\cite{J-H} and  Section 4.11 of \cite{Reu}).
It is a selfadjoint extension of the
following operator, initially defined  in a dense subspace of
${\cal H}_{ph} \otimes {\cal H}_{sp} $,
\be\label{7.1} H(h) = H_{ph} \otimes I + h H_{int},\ee
where $H_{ph} =   h {\rm d} \Gamma (M_{\omega})$ is the photon free energy operator, acting in
a domain   $D(H_{ph}) \subset {\cal H}_{ph}$ and with
\be\label{7.8} H_{int} = \sum _{\lambda =1}^N  \sum _{m=1}^3 (\beta_m + B_m (x_{\lambda})) \otimes
\sigma_m^{[\lambda]},\ee
where  $\beta = (\beta _1 , \beta _2 , \beta_3)$  is the external constant magnetic
field and the  $x_{\lambda }$ ($1\leq \lambda \leq N$)  are the points of
$\R^3$ where the fixed particles are located.

  If an element $U$ of $H^2$  lies in the domain  $ D( M_{\omega }^{-1/2})$
then the Segal  field  $\Phi_S(U)$ is bounded from  $D(H_{ph})$ into ${\cal H}_{ph}$, see point ii) of Proposition 3.4 in
\cite{A-L-N-2} or see {\cite{DG}}.
This is therefore the case for the operators $B_j (x)$ and $E_j(x)$ according
to the assumptions on the ultraviolet cut-off function $\chi$ in  (\ref{7.3}).
Thus, according to Kato-Rellich Theorem,   $H(h)$ has a selfadjoint
extension with the same domain as the free operator $H_0 = H_{ph } \otimes I$
domain.

The Wick symbol definition of some operators is now recalled.

\begin{defi}\label{d-EC-Wick}   Let $H$ be any arbitrary infinite dimensional separable Hilbert space. The coherent states $\Psi_{X , h}$  are the elements of the Fock space ${\cal F}_s (H_{\bf C})$ defined for 
any  $X= (a , b)  $ in $H^2$ and every $h>0$ by,
     \be\label{1.2} \Psi_{X h} = e^{-{|X|^2  \over 4h}} \sum _{n\geq 0} \frac
     {(a+ib) \otimes \cdots \otimes (a+ib)}  {(2h)^{n/2} \sqrt {n!} } \ee
 (with $n$ tensor products in the above sum). When $X=0$, $\Psi_{0 , h} $ denotes the (quantum) vacuum state.
 
 If $A$ denotes any bounded operator in ${\cal F}_s (H_{\bf C})$ or any unbounded
operator $(A , D(A))$ with a domain $D(A)$ containing the coherent states then
the Wick symbol of $A$ is the  function  $\sigma_h^{wick} (A)$ defined on
$H^2$ by,
\be\label{1.3}\sigma_h^{wick} (A) (X) = < A \Psi_{X , h} ,    \Psi_{X , h}>,
\ee
where, for all  $X= (a , b)$ in  $H^2$ and every $h>0$, $\Psi_{X , h}$ is the
corresponding coherent state recalled above. It is also called Husimi function.

 If $A$ is a similar operator in ${\cal F}_s (H_{\bf C})
\otimes {\cal H}_{sp} $  where ${\cal H}_{sp} $ is a finite dimensional Hilbert
space, then the Wick symbol of $A$ is the mapping  $\sigma_h^{wick} (A)$ defined from
$H^2$ into ${\cal L} ( {\cal H}_{sp} )$ verifying, 
\be\label{1.3-bis} < \sigma_h^{wick} (A) (X)a, b>  = < A (\Psi_{X , h}\otimes a)  ,
   \Psi_{X , h} \otimes b>,
\ee
for all $a$ and  $b$
in ${\cal H}_{sp}$.
\end{defi}

For Wick symbols, see {\it e.g.}, \cite{F}, \cite{CR} and also
\cite{Hepp}, \cite{A-N}, \cite{AKN}, $\dots$.  Also recall that the parameter $h$ is the semiclassical
parameter.

If $A$ is a Segal field, its Wick symbol is well defined since coherent states belong to the number operator domain being itself included in the domain of any Segal field (see Lemma 2.1 in  \cite{DG}). This also holds true for the operators ${\rm d} \Gamma (T)$ where $T\in {\cal L} ({\cal H})$.  The Wick symbols of these operators are considered in Theorem \ref{quantif-Wick}.
One denotes by  $X$ or $(q, p)$ the variable of $H^2$ and by $x,k$  variables of
$\R^3$. The Wick symbol of an operator $A$ will be here often denoted by
$A(X, h)$ or $A(q, p, h)$.

The aim of this work consists of studying, for some bounded or unbounded observable operators $A$ in ${\cal H}_{ph} \otimes {\cal H}_{sp} $, the average value of the observable $A$ at time $t$ when starting with an initial  state $ \Psi_{X h} \otimes a$ where
$\Psi_{X h} $ is the coherent state associated with $X\in H^2$ and the element 
$a$ belongs to ${\cal H}_{sp}$ with unit norm. This average value can therefore be written as:
$$ < A e^{-i{t \over h} H(h)} ( \Psi_{X h} \otimes a ),
e^{-i{t \over h} H(h)} ( \Psi_{X h} \otimes a ) >.$$
This average value makes sense at least if $A$ is a bounded operator.
This average value is thus read as a Wick symbol. For our purpose, the initial observables $A$ are assumed to have the particular following form:
\be\label{forme-A} A =  \Phi_S(F_A)  \otimes I + I \otimes S_A \ee
where $S_A$ belongs to ${\cal L} ({\cal H}_{sp})$ and  $\Phi_S(F_A)$ is the Segal field associated
with  an element $F_A$ in the domain $ D( M_{\omega }^{-1/2}) \subset H^2$. In this case,
according to the above remarks, the following product
\be\label{A(t, h)} A (t, h) = e^{i{t \over h} H(h)} A e^{-i{t \over h} H(h)} \ee
is well defined as an operator from $D(H(h))$ into
${\cal H}_{ph} \otimes {\cal H}_{sp} $.
According to Proposition \ref{p-classe}, if $A$ has the form (\ref{forme-A}) then the operator
$A (t, h)$ is the sum of a Segal field  with a bounded operator.
Since both of them have a Wick symbol then
 the Wick  symbol $A(t, X, h)$ taking values in ${\cal L} ( {\cal H}_{sp})$ 
 is well defined. The above mean value is $< A(t, X, h) a, a >$.

For observables being as in (\ref{forme-A}) type, we shall for the function  $A(t, X , h)$,
 on the one hand study bounds of the derivatives and on the other hand give an asymptotic expansion as $h$ goes to $0$ with the aim of obtaining quantum corrections for the Bloch model.
In order to obtain bounds on the  derivatives, we define a class of functions 
$F\in C^{\infty}$ on the phase space $H^2$  associated with a  nonnegative quadratic form
$Q$  on $H^2$ in the following way.

\begin{defi}\label{d1.1}   For any
real separable Hilbert space $H$ and for each  nonnegative quadratic form
$Q$  on
$H^2$, let  $S (H^2, Q) $ be the class of functions $f \in C^{\infty }
(H^2)$  such that  there exists $C(f) >0$ satisfying, for any integer $m\geq
0$, for all vectors $V_1$ ... $V_m$ in $H^2$:
\be\label{1.1}|(d^m f ) (x) (V_1 , \dots , V_m ) | \leq C (f)
Q( V_1) ^{1/2} \dots  Q( V_m) ^{1/2}.
\ee
The smallest constant $C(f)$ satisfying
$(\ref{1.1})$ is denoted by $\Vert f \Vert _Q$.

\end{defi}

In what follows, the quadratic form $Q$ will be,
\be\label{1.100}Q(X) = ( A_QX ) \cdot X , \ee
with $A_Q\in {\cal L}(H^2)$, selfadjoint,   nonnegative, trace class in
$H^2$. Note that the idea of defining a class of symbols  this way for 
quantization purposes, thanks to a quadratic form on the phase space, goes back to
H\"ormander \cite{HO} and Unterberger \cite{U}. 
One however notes that, concerning these works, the constant $C$ denoted by $C_m$ involved there depends on  $m$ whereas it is independent of $m$ here.
Also observe  that a function in  $S(H^2, Q)$ extends as a  holomorphic function
on the complexified  $(H_{\bf C})^2$, satisfying
$$ |F(Z)|\leq \Vert F \Vert _Q \ e^{ Q( {\rm Im} Z) ^{1/2}},\quad
Z\in (H_{\bf C})^2.$$

For each time $t$, we shall show that the function
 $X \rightarrow A(X , t, h)$ belongs to the class of Definition \ref{d1.1} associated with the time dependent quadratic form $Q_t$ on $H^2$ that we now define.

We remark that,
\be\label{HphEC}e^{-i{t\over h} H_{ph} }  = \Gamma (\chi_t) \ee
where $\chi_t$ is the symplectic map defined, setting  $\omega (k) = |k|$, by
\be\label{def-chi-t} \chi_t (q , p) = (q_t , p_t), \ee
$$q_t (k) = \cos (t \omega (k)) q(k) +  \sin (t \omega (k)) p(k),\qquad
p_t (k) =-  \sin (t \omega (k)) q(k) +  \cos (t \omega (k)) p(k).$$


For each operator $A$, bounded from $D( H_{ph}\otimes I)$ to
${\cal H}_{ph} \otimes {\cal H}_{sp}$, we set:
\be\label{evol-libre}
A^{free } ( t, h) =   e^{i{t \over h}(H_{ph} \otimes I)}
A
e^{- i{t \over h} (H_{ph} \otimes I) }.
\ee
In particular, we define in this way $H_{int} ^{free } (t)$.
The proof of Theorem \ref{t-Bogoliub} shows that the Wick symbol
$H^{free}_{int}(t , X)$ of the operator   $H_{int} ^{free } (t)$
is well defined and that:
\be\label{wick-free} H^{free}_{int}(t , X) = H_{int}( \chi_t (X))
. \ee
For all  $t\in \R$, one defines a quadratic form  $Q_t$  on $H^2$ by:
\be\label{def_Q_t} Q_t (q, p) =  |t|  \int_0^t
| dH^{free}_{int}(s , q, p) |^2
ds = |t|  \int_0^t
| dH_{int}(\chi_s( q, p) ) |^2
ds \ee
where the norm in the integral is that of  ${\cal L} ( {\cal H}_{sp})$. For this equality to
define a  quadratic form, one chooses on ${\cal L} ( {\cal H}_{sp})$  the
Hilbert-Schmidt norm. One denotes by  $dH^{free}_{int}(t , q, p)$
the differential with respect to $(q, p)$ of this  affine  function, that is to say the function
obtained by replacing the constant field $\beta$ by $0$.
Notice that the operator  $A_t$ satisfying  $Q_t (X) =
(A_tX)\cdot X$ for all  $X$ in  $H^2$ (as in (\ref{1.100}))
is trace class.

Now we can state our main result.

\begin{theo}\label{tPrincipal} Let $A$ be an observable of the form
(\ref{forme-A}) with $F_A$ in the domain $ D( M_{\omega }^{-1/2})$.
Let $A(t, h)$ be the operator defined in (\ref{A(t, h)}) and $A(t, X , h)$ be
its Wick symbol. Then,

i) The function $X\rightarrow A(t , X, h)$ is the sum of a linear function of the variable $X$ with a function belonging to $S(H^2 , 4 Q_t)$ and taking values in ${\cal L} ({\cal H}_{sp}) $.

More precisely, there exists a sequence of functions 
 $A^{[j]} (t,  \,\cdot\,)$,
$j\geq 0$,  taking values in  ${\cal L} ({\cal H}_{sp}) $ and  satisfying the following properties:

ii) The function $A^{[0]} (t, \,\cdot\,)$ is the sum of a linear function with a function lying in $S(H^2 , 4 Q_t)$ and taking values in ${\cal L} ({\cal H}_{sp}) $ with a
norm bounded independently of $t$
and $h$ when $t$ remains in a compact subset of $\R$ and $h$
varies in  $(0,1]$.

iii) For $j\geq 1$, the function  $A^{[j]} (t, \,\cdot\,)$  belongs to  $S(H^2, 16^{j+1} Q_t)$
and takes values in ${\cal L} ({\cal H}_{sp})$ with a norm  bounded independently of $t$
and $h$  when $t$ remains in a compact subset of $\R$ and $h$
varies in  $(0,1]$.

iv) The Wick symbol   $X \rightarrow A(X , t, h)$ satisfies for any integer  $M$:
\be\label{expansion} A(X , t, h)  = \sum _{j= 0}^M  h^j A^{[j]} (t, X)+
h^{M+1} R^{[M]} (t, X, h) \ee
where $ R^{[M]} (t,  \,\cdot\,, h)$ belongs to $S(H^2, 16^{M+5} Q_t)$ with a norm
bounded independently of $t$
and $h$ when $t$ remains in a compact subset of $\R$ and $h$
varies in  $(0,1]$.

\end{theo}

The construction of the functions  $A^{[j]} (t, \,\cdot\,)$ for $j\geq 0$ and $ R^{[M]} (t, \,\cdot\,, h)$ for any integer $M$ and $h>0$ is given in Section \ref{s-5}.

We can be more precise about  the expansion when the observable $A$ of Theorem
\ref{tPrincipal}  is one of the Segal fields   $B_m(x) = \Phi_{S, h} (B_{mx})$
or $E_m (x) = \Phi_{S, h}  (E_{mx})$,
or one of the operators $\sigma_m ^{[\lambda]}$
defined in (\ref{spin-ini}). These operators all have  the form  (\ref{forme-A}).
We then denote by $B_m^{[j]} (x, t, X)$, $E_m^{[j]} (x, t, X)$ and $  S_m^{[\lambda, j]}  ( t,X)$
($1 \leq m \leq 3$)  the functions denoted by     $A^{[j]}(X , t) $  in
Theorem \ref{tPrincipal}. The detailed construction of these functions  is given in Theorem \ref{tCalcul-termes}. Let us now give
 only the general idea.

The first terms $B_m ^{[0]} (x, t, X)$ and $E_m ^{[0]} (x, t, X)$
are  functions of $X\in H^2$ and also of $(x, t)$. As functions of $(x , t)$, these
functions satisfy the free Maxwell equations,  the initial conditions being the fields
corresponding to the initial coherent
state.  Then, the first terms $  S_m^{[\lambda, 0]}  ( t,X)$
satisfy the Bloch equations \cite{BL} (1946) but where the magnetic field is the sum
of the constant external field with ${\bf B} ^{[0]}$.

Next, the terms with $j\geq 1$ are determined by induction. Let us assume that  
$B_m^{[j-1]} (x, t, X)$ and $E_m^{[j-1]} (x, t, X)$ together with  $  S_m^{[\lambda, j-1]}  ( t, X)$
are already determined. Then, the functions $B_m^{[j]} (x, t, X)$ and $E_m^{[j]} (x, t, X)$ considered as functions of  $(x , t)$ satisfy  Maxwell equations with entirely vanishing initial
condition together with a zero charge density and a current density equal to
$$ {\bf J}^{[j]} (x , t, X) =  \sum _{\lambda = 1} ^N {\bf S^{[\lambda , j-1]}}
(t, X) \times {\rm grad }
\rho (x-x_{\lambda})$$
with
\be\label{rho}
\rho (x) = (2 \pi )^{-3}  \int _{\R^3} |\chi (|k|)|^2  \ \cos (k \cdot  x) dk,\ee
where $\chi$  is the ultraviolet cut-off  function appearing in Definition \ref{7.3}
of the magnetic field operators. This term expresses the radiation, between
times $0$ and $t$, of the spins in the $(j-1)$-th order of their
movement.  One finally determines  $S_m^{[\lambda, j]}  ( t, X)$
solving the differential system (\ref{MaxSym-6}).
It depends on the one hand, on the
mutual interactions of the spins and on the other hand, on quantum corrections
coming from QED. This constitutes a quantum correction of the Bloch equations.
The proof relies on some equations  of Maxwell-Bloch type  for operator valued functions
(see  \cite{SP} and also Theorem \ref{tSpohn}).

We now turn to the number operator $N$ time evolution. The number operator is not under the form 
 (\ref{forme-A}) and one cannot therefore directly apply Theorem \ref{tPrincipal}.
Moreover, the operator
$$  N(t, h) =    e^{i{t \over h}   H(h) } ( N \otimes I)  e^{-i{t \over h} H(h)}. $$
is not precisely defined. We shall instead use its formal derivative:
 \be\label{N'(t)} N'(t, h) =  (i/ h)    e^{i{t \over h}   H(h) }
 [ H(h) , ( N \otimes I) ]  e^{-i{t \over h} H(h)}. \ee
This definition makes sense.
We shall see  that
$[H(h) , N\otimes I ]$ is a bounded operator from $D( H(h))$ into
${\cal F}_s ( H_{\bf C}) \otimes  {\cal H}_{sp}$.
In particular, the quantity $< N'( t, X,  h) a , a>$ amounts to the photon number average value emitted by unit of time, at time $t$, when the initial state is  taken as $\Psi_{Xh} \otimes a$ with a unit normalized element $a$
in  ${\cal H} _{sp}$.

\begin{theo}\label{t-evol-N}

i) The operator $ N'(t, h)$ defined in (\ref{N'(t)}) is well defined from $D( H(h))$ into
${\cal F}_s ( H_{\bf C}) \otimes  {\cal H}_{sp}$. Its Wick symbol $N'(t, X, h)$
defined for $X\in D(M_{\omega} ) \subset H^2$ satisfies,
\be\label{asym-N'-1} N'(t, X, h)=   \sum _{\lambda = 1}^N \sum _{m=1}^3   E_m^{free, pol} (x_{\lambda }, X, t)
  S_m  ^{[\lambda , 0 ]} (t, X )+ N^{res} (t, X, h) ,  \ee
where $X \rightarrow  E_m^{free, pol} (x_{\lambda }, X, t)$ is a linear form on $H^2$
that is determined in Section \ref{s-7} and $X \rightarrow  N^{res} (t, X, h)$ is function in  $S(H^2, K Q_t)$, with $K>0$ and its norm
 in this class is bounded independently of $t$ belonging to a compact set of $\R$
 and of $h$ lying in $(0, 1]$.

ii) There exists a sequence of functions $ N^{[j]} $ on
$\R \times H^2$ satisfying, for  any integer  $M\geq 1$, 
\be\label{asym-N'-2}  N^{res} ( t, X,  h) = \sum _{j=1}^M h^j
 N^{[j]}  (t, X)+ h^{M+1} R^{[M]} (t, X, h).\ee
For all $j\geq 1$, for every $M\geq 1$  and for any $t\in \R$, the function $X \rightarrow N^{[j]} (t, X)$ belongs to $S(H^2, K_j Q_t)$ and
$R_M (t, \cdot , h)$ belongs to a class $S(H^2, L_M  Q_t)$
where  $K_j$ and $L_M$ are some constants and where the norms are bounded independently
 of $t$ belonging to any compact set of $\R$
and of $h$ lying in $(0, 1]$.

\end{theo}

This theorem is proved  in Section \ref{s-7}  where the
functions $N^{[j]}  (t, X)$ are determined.
However, for a better understanding at this stage of equality  (\ref{asym-N'-1}), let us introduce here the element  $ E_m^{free, pol} (x_{\lambda }, X, t)$ involved here. We  define  $ E_m^{free, pol} (x_{\lambda }, X, t)$ in a particular case, when the photon  $X = (q,p)$
satisfies for almost all  $k\in \R^3$,
$$(k\times q(k)  , k\times p(k))  = \varepsilon |k| (-p(k), q(k))$$
where $\varepsilon =  \pm 1$.
These two cases  correspond to the circular right and left polarization
notions in physics. In both cases we have,
$$  E_m^{free, pol} (x_{\lambda }, X, t) = \varepsilon
 {\bf E}^{free}  (x_{\lambda }, t, X). $$
Thus, in other words, (\ref{asym-N'-1}) says that, the mean number of photons emitted by unit
of time is related to the scalar product of the spin and of the electric field, corrected
according to the polarization, up to corrections in ${\cal O}(h)$.

Note that we are not expecting to consider in these issues the limit as t goes to infinity in the semiclassical context recalling that even in finite dimension, semiclassical expansions should be valid up to Ehrenfest time. We also note that \cite{CTDRG} gives first order quantum corrections for a related model, namely, the interaction of an electric dipole moment with the quantized electric field. Besides,  it is observed in \cite{Reu} that the polarization is involved in the  photon emission.  

The relevance of the semiclassical limit is also suggested by the following
 complementary remark. We point out  that the semiclassical approximation of the model studied here not only concerns time evolutions of observables but other issues can naturally be studied. Indeed, it is known, see for example \cite{A-H}\cite{G}\cite{H-S} (note that very often $h=1$ in references for the third model), that the Hamiltonian
 $H(h)$ defined above has a ground state $U_h$ satisfying:
  $$ H(h) U_h = E_h U_h,\quad  E_h = \inf \sigma (H(h)),\quad \Vert U_h \Vert = 1 $$
and it is proved in \cite{A-N-CS} under some conditions that the observable average values
  $B_m(x)$ and $E_m (x)$ related to the three components of the magnetic and electric fields at an arbitrary point  $x\in\R^3$ satisfy,
  $$ < B_m(x) U_h , U_h > = h  B_m^{class } (x) + {\cal O}(h^{3/2}),\quad
< E_m(x) U_h , U_h > = h  E_m^{class } (x) + {\cal O}(h^{3/2})$$
where ${\bf B }^{class } (x) $ and ${\bf E }^{class } (x) $ are the magnetic and electric fields associated, according to elementary physics, with the  $N$ spins viewed as magnets all pointing in the direction of the non quantized constant external magnetic field. In particular,
 ${\bf E }^{class } (x)= 0 $. We also derive in \cite{A-N-CS} a connection between the model studied here (the third model) with the Ising model. 
 
 {\it Notations.} The scalar product of two elements $a$ and $b$ of a real Hilbert space will be
denoted by  $a \cdot b$. It is the case for
the configuration space $H$ or the phase space $H^2$, when it is not identified
with the complexified $H_{\bf C}$. We denote  by  $<a , b>$ the  hermitian
product on a complex space, which always will be antilinear with respect to the second factor
(such that $< a , ib>= -i <a, b>$).   It is the case for  $H_{\bf C}$, for the
Fock space ${\cal H}_{ph} = {\cal F}_s (H_{\bf C})$, and for the space ${\cal H}_{sp}$.
The $H$ or $H^2$ norm will be denoted by  $| \cdot |$ and the
norm of  ${\cal F}_s (H_{\bf C})$  or $L^2 (B , \mu _{B , h/2})$, by
$\Vert \cdot \Vert $.

{\it Sketch of the proof.}  
One first writes the equations that should be satisfied by the functions   $A^{[j]} (t, X)$ in Theorem  \ref{tPrincipal} in order to be, at least formally, a good approximation of the Wick symbol of the operator   
$A(t, h)$.
These equations are explicitly written down in Sections \ref{s-A0} and \ref{s-Aj},
where we also prove existence and unicity properties of the  solutions. It is in addition derived that these solutions belong to some classes of Definition \ref{d1.1}. It remains to prove 
(\ref{expansion}) and thus, to compare the true Wick symbol  of $A(t, h)$ with its supposed approximation. Rather than comparing two functions, it seems easier to compare two operators.
Therefore, with each function $F$ belonging to a class  $S(H^2, Q)$ of Definition
\ref{d1.1}, one associates an operator denoted  by $Op_h^{wick} (F)$ whose Wick symbol is $F$.

This Wick quantization of a given function $F$ is not always possible. It is possible if $F$ is polynomial function. The corresponding operator is then defined using Wick (normal) ordering (see \cite{A-N}, \cite{G-J}, $\dots$). It is also defined for some quadratic forms (see \cite{LA}).

We show in Section \ref{s-2} that Wick quantization is also possible for functions $F$ in a class  $S(H^2, Q)$ of Definition \ref{d1.1}.  
To this end, we first begin by giving an heat inverse operator to the function $F$.
The fact that this is
possible should probably be related to  properties of the functions in $S(H^2 , Q)$,
which are stronger than analyticity. Next, we  use the anti-Wick quantization commonly used in finite dimension (see \cite{B-S} and \cite{Ler}) which is however only an intermediate step in our case.
Also in infinite dimension, this anti-Wick quantization has in own interest 
 (see  \cite{A-N} for constructing semiclassical measures).
This will enables a Wick quantization of the coefficients 
$X\rightarrow  A^{[j]} (t, X)$ and also the error terms appearing in computations.  Thus, estimate (\ref{expansion}) brings us back to a comparaison between two operators which is easier to consider. Let us also mention that another technique of \cite{K-R} concerning Wick quantization  could may be applied to our classes.

\section{Quantization}\label{s-2}

The purpose of this section is to give an answer to the following issue when $H$ is any   arbitrary  infinite dimensional  separable Hilbert space not necessarily being the photon Hilbert spaces.
For a given function $F$ on $H^2$, can we find an operator $A$ in the Fock space ${\cal F}_s (H_{\bf C})$, bounded or unbounded, with a Wick symbol equal to $F$ ?
If so, we can say that $A$ is the Wick quantization of $F$. Before, we had to define an anti-Wick quantization which here is only a first step but can have is own interest.

\subsection{Anti-Wick quantization.}\label{s-AW}

We recall that in finite dimension $n$, the anti-Wick operator associated with a bounded measurable function $F$ on  $\R^n$ is defined by, for all  $f$ and $g$ in
 $L^2(\R^n)$,  for any $h>0$,
\be\label{def-AW-Rn}< Op_h^{AW} ( F)f , g> = (2 \pi h)^{-n}  \int _{\R^{2n} }
F(X)  < f , \Psi _{X , h} > \ <\Psi _{X , h}, g >
dX,\ee
where the  $\Psi _{X , h}$ are the standard coherent states on 
$\R^n$ indexed by $X\in \R^{2n}$ (see \cite{F}). One of the advantages  of this quantization used for example in  \cite{Ler} is the possibility to consider less regular  functions $F$.

With the aim of  defining corresponding operators in the Fock space ${\cal H}_{ph}$, one could imagine to consider integrals on $H^2$. Naturally, the Lebesgue measure is not existing anymore and it can be replaced by a Gaussian measure. For this purpose, $H$ is replaced in the integral by another space $B$.

{\it Wiener Space. }
Let us recall that we can associate with any Hilbert space  $H$ a Banach space $B$ containing
$H$ and, for each $t>0$, a Gaussian measure $\mu _{B , t}$ on the Borel
$\sigma-$algebra of $B$, with the properties of the theorem
below. We say that the  variance is $t$. For these points, see \cite{G-2},
\cite{KU} and also Theorem \ref{t1.1}. The exact conditions to be fulfilled by
$B$ together with the properties involved in this paper are recalled
in \cite{A-J-N}.

\begin{theo}\label{t1.1} (Gross \cite{G-1}-\cite{G-4}, Kuo \cite{KU}).
Let $H$ be a real separable Hilbert space. Then, there exists a (non unique)
Banach space $B$ containing $H$,
such that  $B' \subset H' = H \subset B$, each space being dense in the
subsequent, and for all $t>0$, there exists a probability measure
$\mu_{B , t}$ on the Borel $\sigma$-algebra of  $B$ satisfying,
\be\label{1.1b} \int _B e^{i a(x) } d\mu _{B, t} (x) = e^{-{t\over 2} |a|^2 },\quad a\in B'.\ee
Here, $|a|$ denotes the norm in  $H$ and the notation $a(x)$ stands for
the duality between $B'$ and $B$.
\end{theo}

One says that the triplet $(i, H , B)$ is a Wiener space where $i$ is the
injection from $H$ into $B$. One also says that  $B$ is a Wiener extension of $H$.

{\it Gaussian variables. } We  recall (\cite{KU}) that,  for all $a$ in
$B'\subset H$, the mapping $B \ni  x \mapsto a(x) $ belongs to
$L^2(B , \mu _{B , h})$, with a norm equal to ${h}^\frac{1}{2} |a|$.
Thus, the mapping  associating with every $a\in B'$, the  above function
considered as an element of $L^2(B , \mu _{B , h})$,
can be extended by density to a mapping
$a \mapsto \ell_a$ from $H$ into $L^2(B , \mu _{B , h})$.

{\it Stochastic extensions.} We would like to define an operator starting with some functions $F$ defined on the phase space $H^2$ whereas the integral
formulas  concerns only functions on $B^2$ which is the only  measurable
space available. The stochastic extension of $F$ can be then here implied. The notion of stochastic extension in $L^p (H^2 ,\mu_{B , h})$
of some functions defined on $H^2$ (and satisfying suitable conditions)
is recalled in Definition 2.2 of \cite{A-L-N-3}.
One may find in  \cite{A-L-N-3} and in \cite{A-J-N} examples of functions
admitting such extensions.

{\it  Segal-Bargmann transform.} For each  function $f$ in
${\cal H}_{ph}$, for every $X$ in $H^2$ and for each $h>0$, we set:
\be\label{2.3} (T_{h}f)(X) = \frac { < f , \Psi_{X,h} > } { < \Psi_{0,h} , \Psi_{X,h}> }
= e^{\frac {|X|^2} {4h} } < f , \Psi_{X,h} >  \ee
where $\Psi_{X,h}$ denotes the coherent state defined in  (\ref{1.2}). This
function ${T_h f}$ is G\^ateaux anti-holomorphic on $H^2$ when one identifies
$(x,\xi)\in H^2$  with $x+i\xi$. Its restriction to any finite dimensional
subspace $E$ of $H^2$ belongs to $L^2(E , \mu _{E, h})$, its norm being bounded
independently of  $E$. According to Theorem 8.8 in \cite{A-J-N}, its stochastic
extension $\widetilde T_hf $ exists in $L^2(B^2,\mu _{B^2, h})$ and
$ \widetilde T_h$ is a partial isometry from  ${\cal F}_s(H_{\bf C})$ in
$L^2(B^2,\mu _{B^2, h})$ whose range is the closure of anti-holomorphic functions.
The mapping $\widetilde T_h$ is called the Segal Bargmann transform.

{\it Anti-Wick Operator.} For every function $F$ on $H^2$ (called the symbol) admitting a
stochastic extension  $\widetilde F$ measurable and bounded on $B^2$,
one denotes by  $Op_h^{AW} (F)$ the operator defined by, 
\be\label{def-AW}< Op_h^{AW} ( F)f , g> = \int _{B^2} \widetilde F(X) \widetilde T_hf (X) \overline { \widetilde T_hg (X)}
d\mu _{B^2, h} (X),\ee
for all $f$ and $g$ in
 ${\cal H}_{ph}$.
Since $ \widetilde T_h$ is a partial isometry from  ${\cal F}_{s}(H_{\bf C}) $
into $L^2(B^2 , \mu _{B^2, h})$, this operator is bounded on  ${\cal F}_{s}(H_{\bf C}) $
and its norm is smaller than the sup norm of $\widetilde F$.

Unlike  in the finite dimensional case, we cannot define an anti-Wick
operator for any bounded measurable function $F$ in $H^2$ since this function
needs an stochastic extension. In particular, one can associate an anti-Wick
operator with every function $F$
in $S(H^2, Q)$  where  $Q$ is a nonnegative quadratic form on $H^2$ written under the form  (\ref{1.100}) with $A_Q$ being selfadjoint and trace class.
Indeed, according to \cite{J}, Proposition 3.10, such a function
admits a stochastic extension $\widetilde F$ measurable and bounded on $B^2$. In this case, we
have $\Vert \widetilde F \Vert _{L^{\infty}} \leq \Vert F \Vert _Q$ and therefore, since
$ \widetilde T_h$ is a partial isometry from  ${\cal F}_s(H_{\bf C})$ into
$L^2(B^2,\mu _{B^2, h})$, we have:
\be\label{norme-AW} \Vert Op_h^{AW} ( F) \Vert \leq \Vert F \Vert _Q \ee
Again, we shall use anti-Wick quantization as an intermediate   step towards the Wick quantization studied in the following subsection. The anti-Wick quantization in infinite dimension has its own interest: it is used in  \cite{A-N} for the construction of semiclassical measures.

\subsection{Wick quantization.}\label{s-2-2}

Any bounded operator $T$ has a Wick symbol but it not always possible for a  function $F$ on $H^2$ to find an operator whose Wick symbol is $F$. The  theorem below shows that it is true in some cases. We recall  here that the coherent states belong to all $D(N^m)$ for any integer $m\geq 1$. Polynomial functions are finite linear combinations of multiplications of $e_j\cdot q$ and $e_j\cdot p$ as maps of $(q,p)$.

\begin{theo}\label{quantif-Wick}

 i) For any function $F$ belonging to   $S (H^2, Q) $ where $Q$ is a 
 nonnegative quadratic form on  $H^2$ written under the form  (\ref{1.100}), where $A_Q$ is selfadjoint and trace class,
and for each $h>0$, there exists a unique bounded operator $B_h$ in
${\cal F}_s(H_{\bf C})$ such that $\sigma_h^{wick } (B_h ) = F$. We then set
$B_h = Op_h^{wick } (F)$. On has the following estimate:
\be\label{norme-Wick}  \Vert Op_h^{wick } (F) \Vert \leq e^{ (h/2)   {\rm Tr } A_Q} \Vert F \Vert _Q,  \ee
where the norm in the left hand side is the  ${\cal L} ( {\cal F}_s (H _{\bf C}))$ norm and  $A_Q$ is the operator verifying $Q(X) = (A_Q X) \cdot X$ for all $X\in H^2$.

ii) For any polynomial function  $F$ on $H^2$ of degree $m$ and for all $h>0$ there exists an unbounded operator $B_h$, closable, with initial domain $D(N^ m)$,
satisfying $\sigma_h^{wick } (B_h ) = F$. We thus set
$B_h = Op_h^{wick } (F)$. In particular, the Wick symbol of a Segal field 
$\Phi_S (A)$ associated with $A \in H^2$
 is the function  $H^2 \ni X \rightarrow h^{-1/2} A \cdot X $, where  $A \cdot X $ is the
 real scalar product on $H^2$.

iii) For any quadratic form  $F$ on $H^2$ written as $F(q , p) = (Aq)\cdot q +
(Ap)\cdot p$, where $A\in {\cal L} (H)$ is selfadjoint,  and for every $h>0$,
the operator $B_h= 2h {\rm d} \Gamma (A) $ satisfies  $\sigma_h^{wick } (B_h ) = F$.

\end{theo}

 We  now need to  involve the heat operator in order to first prove point i) of the above theorem.
  The heat operator is defined for each
measurable bounded function $F$ in  $H^2$ admitting a stochastic extension
$\widetilde  F$ in $L^1 (B^2 ,\mu_{B , h})$   and for each $h>0$  by:
\be\label{Ht}(H_{h} F) (X) = \int _{B^2} \widetilde  F (X+ Y)d\mu_{B^2,h} (Y),\ee
for $X\in H^2$ (see Definition 5.1 and formula (28)  of \cite{J}).
We  also have:
\be\label{Ht2} (H_{h} F) (X) =e^{-\frac {|X|^2} {2h} }  \int _{B^2}
\widetilde F (U)\ e^{\frac {1} {h} \ell _X (U)} d\mu_{B^2,h} (U). \ee

\begin{theo}\label{t4.1}  Let  $F$ be in $S(H^2, Q)$. Then the function
$h \rightarrow H_hF$ is $C^{\infty }$ on  $[0, \infty [$  with values in the
Banach space $S (H^2, Q)$. Moreover, for all $h>0$, the operator $H_h$ is
an isomorphism from  $S (H^2, Q)$ onto itself. As operators in  $S (H^2, Q)$, the norm of $H_h$ and the norm of  its inverse
denoted by $H_{-h}$ satisfy the following bounds:
\be\label{4.1} \Vert H_h \Vert \leq 1,\qquad \Vert H_{-h} \Vert \leq
e^{ (h/2)   {\rm Tr } A_Q},
\ee
where $A_Q$ is the operator satisfying $Q(X) = (A_Q X) \cdot X$.

\end{theo}

{\it Proof. } The first claim is Theorem 5.17  of \cite{J}. It is proved in
Proposition 5.12  of \cite{J} that the Laplace operator $\Delta $ is well defined
on $S (H^2, Q)$ and that it is  bounded on
$S (H^2,  Q)$ into itself, with norm smaller or equal than $ {\rm Tr}A_Q$. On the other hand, one
has  $H_h = e^{ (h/2)\Delta}$. This is proved in  \cite{J} (proposition 5.13 or
Theorem 5.17). Therefore, the  inverse of the bounded operator $H_h$  is defined by:
\be\label{4.3} H_{-h} F = \sum _{m=0}^{\infty} (-1)^m \frac { h^m} { 2^m m!}
\Delta ^m F. \ee
Its norm satisfies  (\ref{4.1}), which proves the result.

{\it Proof of point i) of Theorem \ref{quantif-Wick}.}
One knows that, for all $X$ and $Y$ in $H^2$:
\be\label{PSEC} < \Psi_{X h} , \Psi _{Yh}> =e^{-{1\over 4h}(|X-Y|^2 ) + {i\over 2h} \sigma (X , Y) }
\ee
where  $\sigma $  is the symplectic form $ \sigma ( ( x , \xi) , (q, p) ) =
q \cdot \xi - p \cdot x$. Hence, for each $U= (a, b)$ in $H^2$, setting $\check U = (-b , a)$:
\be\label{SBEC} \widetilde T_h \Psi_{U, h} (X) =
e^{ - {1 \over 4h}(|U|^2  + {1 \over 2h}\ell _{ U - i \check U} (X)}.\ee
Then, for each   $F$ in $S(H^2, Q)$ and $h>0$, one has, according to (\ref{def-AW}):
$$< Op_h^{AW} (  F) \Psi_{U, h} ,  \Psi_{U, h} > =
e^{ - {1 \over 2h}|U|^2 }\int _{B^2} \widetilde F (q, p) e^{ {1 \over h}(\ell _a (q)
+ \ell _b (p)) }
d\mu _{B^2 , h} (q, p) =  H_h F(U). $$
The last equality follows from  (\ref{Ht2}).     Therefore:
\be\label{AW-WK} \sigma _h^{wick} (Op_h^{AW} (F)) = H_h F. \ee
To prove point i) of Theorem  \ref{quantif-Wick}, it is then sufficient to set 
\be\label{def-OPW}  Op_h^{W}(F) =   Op_h^{AW} (  H_{-h}  F). \ee
The norm estimate (\ref{norme-Wick}) is then a consequence of the above definition, with  (\ref{4.1}) and (\ref{norme-AW}).

{\it Proof of point ii) of Theorem  \ref{quantif-Wick}. }
If $F$ is a polynomial function of degree $m$ on $H^2$, we can again define $Op_h^{W}(F)$ by (\ref{def-OPW}). Indeed,
  $G=  H_{-h}  F$ is also a polynomial function $G$ on $H^2$ with the same degree. One proves in \cite{J} (Proposition 3.17) 
 that each polynomial function $G$ of degree $m$ has a stochastic extension  $\widetilde G$. 
Moreover,
$$ |(G T_hf) (q , p) | \leq C |(T_h f) (q-ip) | \sum |(q - ip)^{\alpha} |, $$
where the above sum is finite. A Hilbert basis  $(e_j)$ of $H$ is fixed. 
Note that $(q - ip)^{\alpha} (Tf) (q-ip)$ is, up to a multiplicative factor, the Segal Bargmann transform of  $( a^{\star} (e) )^{\alpha } f$, which is well defined as an element of  ${\cal F}_s (H_{\bf C})$   if $f$ belongs to the domain of $N^m$. 
Thus,  the function  $\widetilde G \widetilde T_hf $ belongs to
$L^2(B^2 , \mu _{B^2, H})$. 
Consequently, the integral (\ref{def-AW}) makes sense. It indeed defines an unbounded operator  $ Op_h^{AW} (  H_{-h}  F)$ with (initial) domain $D(N^m)$ where $m$ is the degree of $F$. We recall that the coherent states are in the domain of   $D(N^m)$. Therefore, the Wick symbol of $ Op_h^{AW} (  H_{-h}  F)$
is well defined. Reasoning as in  the above point i) shows that this Wick symbol is equal to $F$.

Let us now see the usual relation between the Wick quantization of a polynomial function with the Wick (normal) ordering.

 {\it Degree one polynomial function.}   First, we prove the last claim of point ii). In order to derive it, we remark, from the Definition
 (\ref{1.2}) of coherent states, that, for each $X\in H^2$:
 \be\label{alter-EC}  \Psi_{X, h} = e^{- \frac {i} {\sqrt {h} } \Phi _S ( \widehat  X )} \Psi_0 \ee
 where $\Psi_0$ is the vacuum (independent of $h$) and $ \widehat X = (-b, a) $
 for $X = (a, b)$. Therefore
 $$ e^{i  t \Phi _S ( A )} \Psi_{X, h} =  e^{i  t \Phi _S ( A )}
  e^{- \frac {i} {\sqrt {h}} \Phi _S ( \widehat X )} \Psi_0. $$
  We know that, for all $U$ and $V$ in $H^2$:
  \be\label{compos-Segal} 
  e^{i   \Phi _S ( U )}  e^{i   \Phi _S ( V )}
  =  e^{ \frac{i}{2} \sigma ( U ,   V  )  }  e^{i   \Phi _S ( U +V )}\ee
 where $\sigma $ is the symplectic form $\sigma ( (a, b) , (q, p)) =
  b\cdot q - a \cdot p$. Therefore:
  $$ e^{i  t \Phi _S ( A )} \Psi_{X, h} = e^{ - \frac {it} {2\sqrt {h}} \sigma ( A ,   \widehat X )  }
    e^{i   \Phi _S ( tA  - (1/ \sqrt {h})  \widehat X )} \Psi_0
    = e^{ -  \frac {it} {2\sqrt {h}}  \sigma ( A , \widehat X )   }
  \ \Psi_{X + t\sqrt {h}\widehat A , h} . $$
According to (\ref{PSEC}),
 $$ < e^{i  t \Phi _S ( A )} \Psi_{X, h} , \Psi_{X, h} > =
 e^{ - \frac {it} {2\sqrt {h}}  \sigma ( A ,   \widehat X )   } \
  e^{-  \frac { t^2} {4} |  \widehat A |^2 } \
   e^{  \frac {i} {2h} {\rm Im} < X +  t \sqrt {h}\widehat A , X > } . $$
 Differentiating with respect to $t$, at $t=0$, and using  $\sigma ( A ,   \widehat X ) = - A \cdot X$ where  $A \cdot X $ is the real scalar product on  $H^2$,   we indeed obtain that:
 $$ < \Phi _S ( A ) \Psi_{X, h} , \Psi_{X, h} > =
 h^{-1/2} A \cdot X. $$
If $F$ is a degree one polynomial function on $H^2$ and if $(e_j)$
is a Hilbert basis  of  $H$, we can write:
$$ F(q, p) = \sum _j a_j (e_j \cdot q + i e_j \cdot p) +
b_j (e_j \cdot q - i e_j \cdot p).$$
Here the sums are finite. Form the above computations:
\be\label{m=1} Op_h^{wick} (F) = \sqrt {2h}  \sum _j a_j  a(e_j) + b_j  a^{\star} (e_j). \ee

{\it  Wick normal ordering. Degree $m$ polynomial functions.} Any polynomial function $F$ of degree $m$ can be written as:
$$ F(q , p) = \sum _{|\alpha | + |\beta | \leq m} a_{\alpha \beta } (  q + ip)^{\alpha }
(  q - ip)^{\beta } $$
setting, for any multi-index  $\alpha $,
$$ (  q \pm  ip)^{\alpha }  = \prod _j ( e_j \cdot q \pm  i e_j \cdot p) ^{\alpha _j}. $$
Let us check that,
\be\label{ordre-wick} Op_h^{wick} (F) = (2h) ^{m/2}  \sum _{|\alpha | + |\beta | \leq m} a_{\alpha \beta }
 a^{\star} (e ) ^{\beta}  a (e ) ^{\alpha } \ee
where we set:
$$  a (e ) ^{\alpha } = \prod _j a(e_j)^{\alpha _j}.$$
This is proved by iteration on $m$. For $m=1$, it is (\ref{m=1}). Assume that the equality is proved for all polynomial function of degree $\leq m-1$. Each degree $m$ polynomial function can be expressed as:
$$ F(q, p) = (e_j. q + i e_j p ) A(q, p) + (e_k. q - i e_k p ) B(q , p) $$
where $A$  and $B$ are of degree $\leq m-1$. Then, the operator
$$ T_h =  \sqrt {2h} \Big ( a^{\star } (e_k) Op_h^{wick } (B ) +
 Op_h^{wick } (A )  a (e_j)\Big )$$
 has $F$ as Wick symbol. This follows from  (\ref{m=1}) and from an analog of  Theorem
\ref{t4-1} valid for finite sums for the Wick symbol of the composition of two  operators. Using the induction hypothesis for $A$ and $B$, we indeed deduce that 
 (\ref{ordre-wick}) holds true for all $m$.

{\it Proof of point iii) of Theorem \ref{quantif-Wick}.}
 By Definition (\ref{1.3}) of the Wick symbol and according to Definition (\ref{1.2}) of the
   coherent states together with the usual definition of ${\rm d} \Gamma (A)$,
    it follows that the Wick symbol of an operator ${\rm d} \Gamma (A)$ where $A\in {\cal L} (H)$ and selfadjoint is given by
 \be\label{SW-dGamma}   \sigma _h ^{wick} ( {\rm d} \Gamma (A) ) (q, p) = {1\over 2h}
 \Big [( Aq ) \cdot  q + ( Ap ) \cdot  p \Big ].\ee

 For instance, the Wick symbol of the number operator  $N= {\rm d} \Gamma (I)$ is:
\be\label{Sym-W-N} N(q , p) ={1\over 2h} ( |q|^2 + |p|^2). \ee
 The Wick symbol of    $H_{ph} = h {\rm d} \Gamma (M_{\omega})$ is defined only for
$X= (q, p)$ in $D(M_{\omega})$, by:
 \be\label{Sym-W-Hph} H_{ph}(q , p) ={1\over 2} \Big [  ( M_{\omega} q ) \cdot q  + ( M_{\omega} p ) \cdot p
 \Big ].  \ee

\section{Beals characterization. }\label{s-3}

Wick quantization implements a one to one correspondance between the space $S(H^2 , Q)$ in Definition \ref{d1.1} and a space of operators acting in  ${\cal F}_{s} (H _{\bf C})$ that we now specify.

This set of operators is inspired by the Beals characterization  \cite{Beals}.
In the finite dimensional case, the multiplication operators  by a coordinate function
and the partial differentiation operators play an essential part. In the
infinite dimensional case, this role is played by Segal fields. The usual Beals hypothesis is here modified in order to avoid 
domain definition issues. 
\begin{defi}\label{d1.3} (\cite{A-L-N-3}).
For each $V= (a, b)$  in $H^2$, let
 $ \widehat V = (-b, a)$. Set $h>0$.
Let  $Q$ be a quadratic form on  $H^2$ of the form (\ref{1.100})
with $A_Q\in {\cal L}(H^2)$   selfadjoint,   nonnegative and trace class in
$H^2$. We denote by  ${\cal L} (Q) $  the space of all bounded operators  $A $
  in ${\cal F}_{s} (H _{\bf C})$  satisfying:

i) For each integer  $m\geq 1$,   for all  $V_1,\dots,V_m$ in $H^2$  and for each
$h>0$, the function
\be\label{3-1} \R^m \ni t \rightarrow   A_h(t) = e^{ \frac {i} {\sqrt {h}} \Phi_S ( t_1  \widehat  V_1+ ...
+ t_m  \widehat  V_m) } \ A \ e^{ - \frac {i} {\sqrt {h}} \Phi_S ( t_1  \widehat  V_1+ ...
+ t_m  \widehat  V_m) } \ee
is  class $C^m$ from $\R^m$ into ${\cal L } ( {\cal F}_{s} (H _{\bf C}))$.

ii) There exists a positive real number $C(A,Q)$ satisfying, for each integer  $m$ and for all  $V_1,\dots,V_m$ in
$H^2$, for any $h$ in $(0, 1]$:
\be\label{1.5} \Vert \partial _{t_1} .... \partial _{t_m} A_h(0)  \Vert  \leq
C(A,Q)    \prod _{j= 1}^m Q( V_j ) ^{1/2},\quad h\in (0, 1].\ee
The norm in the above right hand side is the  ${\cal L } ( {\cal F}_{s} (H _{\bf C}))$ norm.
The smallest constant $C(A,Q) $ such that (\ref{1.5}) is valid is denoted by
$\Vert A \Vert _{ {\cal L} (Q)}$.
\end{defi}

In the usual Beals hypothesis, the estimate (\ref{1.5}) would be written as:
$$ \Vert {\rm ad}   \Phi_S( \widehat V_1) \cdots {\rm ad}  \Phi_S( \widehat V_m) A  \Vert  \leq
C(A,Q) \ h^{m/2}   \prod _{j= 1}^m Q( V_j ) ^{1/2},\quad h\in (0, 1].
$$

If $A$ ad $B$ are in  ${\cal L} (Q) $ then the composition $A \circ B$ belongs to ${\cal L} (4Q) $ and
\be\label{comp}  \Vert A \circ B \Vert _{{\cal L} (4Q)} \leq
\Vert A  \Vert _{{\cal L} (Q)}\
\Vert B \Vert _{{\cal L} (Q)}. \ee

\begin{theo}\label{t.Beals}  Let $Q$  be the quadratic form  written as
 in (\ref{1.100}), with $A_Q\in {\cal L}(H^2)$, selfadjoint,   nonnegative, trace class in
$H^2$. Then:

i) For each  operator $A$ in ${\cal L} (Q)$ and for each $h$ in $(0, 1]$, the Wick symbol
$\sigma _h ^{wick} (A)$ belongs to $S(H^2 , Q)$ and we also have:
$$ \Vert \sigma _h ^{wick} (A) \Vert _Q  \leq \Vert A \Vert _{{\cal L} (Q)}.$$
ii) For each function  $F$  being in $S(H^2 , Q)$ and for every $h>0$, the operator
$Op_h^{wick} (F)$ belongs to ${\cal L} (Q)$ and the following estimates holds true:
\be\label{Beals-3}
\Vert Op_h^{wick} (F) \Vert_{{\cal L} (Q)}   \leq
\Vert F \Vert _Q e^{ (h/2) {\rm Tr } (A_Q)} .\ee

\end{theo}

Theorem \ref{t.Beals} is an infinite dimensional type of  Beals characterization Theorem (see \cite{Beals}, and also \cite{A-L-N-3}).

{\it  Proof of Theorem \ref{t.Beals}. }

i) By (\ref{alter-EC}) and (\ref{compos-Segal}), we have, for all $X$ in $H^2$
and for each finite system $(V_1 , ... V_m)$ in $H^2$:
$$ \sigma _h ^{wick} (A_h(t) ) (X)= \sigma _h ^{wick} (A) (X + t \cdot V). $$
Therefore:
$$ ( d^m \sigma _h ^{wick} (A) ) (X) (V_1 , ... , V_m) =
\sigma _h ^{wick} \Big ( \partial _{t_1} .... \partial _{t_m} A_h(0)  \Big )(X). $$
Point i) then follows since $|\sigma_h^{wick} (B) (X)| \leq \Vert B \Vert$  for any bounded operator $B$ and for all $X\in H^2$.

ii) Let  $F$ be in   $S(H^2 , Q)$,  $h>0$ and  $A= Op_h^{wick} (F)$.
The above computations show that $A_h(t) =  Op_h^{wick} (F ( \cdot + t \cdot V) )$ for any finite sequence
$(V_1 , ... , V_m)$  in $H^2$.
The function $\R^m \ni t \rightarrow F ( \cdot + t \cdot V) $ belongs to the class $C^m $
from  $\R^m$ into $S(H^2 , Q)$. The mapping $Op_h^{wick}$ is a continuous linear map from
$S(H^2 , Q)$ into ${\cal L } ( {\cal F}_{s} (H _{\bf C}))$ with a norm bound given in 
(\ref{norme-Wick}). Consequently, the function $t \rightarrow A_h (t)$ is in the class $C^m$  from $\R^m$ into ${\cal L } ( {\cal F}_{s} (H _{\bf C}))$.  We have:
$$ \partial _{t_1} .... \partial _{t_m} F ( X + t \cdot V)  \Bigg |_{t= 0} =
(d^m F) ( X) ( V_1 , ... V_m).$$
Point ii) then follows.

\section{Compositions, commutators and covariance. }\label{s-4}

\subsection{Mizrahi series.}\label{s-4-1}

The operator composition will be used in two cases. In Section \ref{s-5}, the commutator of two operators is studied, one being a Segal field and the other one being a bounded operator with a Wick symbol belonging to the class  $S(H^2, Q)$. In Section \ref{s-6}, the composition of two bounded operators is considered, when the Wick symbols are in a class $S(H^2, Q)$.

\begin{theo}\label{t4-1}  Let $G$ be a function in  $S(H^2 , Q)$ with a
quadratic form $Q$ with $A_Q$ being nonnegative and trace class.  Let $V$ be an element of $H$.
Then we have:
\be\label{Mizrahi-a} [a(V), Op_h^{wick } (G) ] = \sqrt {h/2}\  Op_h^{wick } (K),\qquad
K(q, p) = ( V\cdot \partial_q + i V\cdot \partial_p) G(q, p). \ee
One has:
$$ \sigma _h^{wick}\Big ( \Phi_{S, h} (V) \circ  Op_h^{wick } (G) \Big )  (X)  =
\varphi_V (X) G(X) + h C^{1, wick } (\varphi_V, G)(X)  $$
where $\varphi_V(X) = V\cdot X$. Fixing a Hilbert basis   $(e_j)$ of $H$,
one sets $\partial _{q_j} = e_j\cdot \partial _{q}$  and for all differentiable functions  $F$ and $G$,
$$ C^{1, wick } (F, G) = \frac {1} {2}
 \sum _{j}
 (\partial _{q_j} - i \partial _{p_j} )  F  \
( \partial _{q_j} + i \partial _{p_j } ) G.  $$
The above sum is independent of the chosen basis. We also have:
$$ \sigma _h^{wick} \Big ( Op_h^{wick } (G) \circ  \Phi_{S, h} (V) \Big ) (X)  =
\varphi_V (X) G(X) + h C^{1, wick } (G , \varphi_V)(X).  $$

\end{theo}

See \cite{A-N}. Fix a Hilbert basis  
 $(e_j)$ of $H$.
We define, for all multi-index  $\alpha = ( \alpha _j) $
(which means $\alpha _j = 0$ except for a finite number of values of $j$),
two differential operators  on $H^2$, denoting by
$(q , p)$ the variable of $H^2$:
$$ (\partial _q \pm i \partial _{p})^{\alpha }
=  \prod _j  \left ( e_j \cdot  \partial_ q \pm i
 e_j \cdot   \partial   _ p\right )^{\alpha _j} . $$

\begin{theo}\label{t-Mizrahi}  Let $Q$ be a quadratic form on $H^2$
  where $A_Q$ is nonnegative and trace class. Then, for each
  $F$ and $G$  in  $S(H^2 , Q)$, we can write:
 \be\label{Miz-1} Op_h^{wick} (F) \circ  Op_h^{wick} (G) =  Op_h^{wick} (C_h^{wick} (F , G)  ) \ee
 \be\label{Miz-2} C_h^{wick} (F , G)  = \sum _{\alpha}
  \left (\frac  {h^{|\alpha | }}{ 2^{|\alpha | }  \alpha !} \right )
(\partial _q - i \partial _{p})^{\alpha }  F  \
( \partial _q + i \partial _{p })^{\alpha } G.
 \ee
The above series is absolutely converging. For all integers $m$, we have:
\be\label{majo-reste} \left | C_h^{wick} (F , G)(X) - \sum _{|\alpha| \leq m }
  \left (\frac  {h^{|\alpha | }}{ 2^{|\alpha | }  \alpha !} \right )
(\partial _q - i \partial _{p})^{\alpha }  F (X) \
( \partial _q + i \partial _{p })^{\alpha } G(X) \right  |  \leq \ee
$$ \Vert F \Vert _{Q }\
\Vert G \Vert _{Q } \  h^{m+1}  [ {\rm Tr } A_Q ]^{m+1} e^{ h {\rm Tr } A_Q }.  $$

\end{theo}

See \cite{Appleby} or \cite{Mizrahi} in the case of finite dimension.
See  also \cite{A-N}   formula (15)(i).

{\it Proof.} For all multi-indices $\alpha$, we define an element $u_{\alpha}$ of
${\cal F}_s (H_{\bf C})$ by:
  \be\label{u-alpha} u_{\alpha} =  (\alpha !)^{-1/2} \left [\prod _j
 \Big ( a^{\star} (e_j) \Big ) ^{\alpha _j} \right ] \Psi_0
 \ee
where $\Psi_0$ is the vacuum state. Thus, $(  u_{\alpha})$ is a Hilbert basis  of ${\cal F}_s (H_{\bf C})$.  We see that:
 \be\label{prod-scal-base} <  Op_h^{wick} (G) \Psi_0 , u_{\alpha} > =
  \left (\frac  {h^{|\alpha | }}{ 2^{|\alpha | }  \alpha !} \right )^{1/2}
( \partial _q + i \partial _{p })^{\alpha } G (0).\ee
Indeed,
$$<  Op_h^{wick} (G) \Psi_0 , u_{\alpha} > = (\alpha !)^{-1/2}
<  \prod _j \Big (  {\rm ad}\  a(e_j)\Big ) ^{\alpha _j}  Op_h^{wick} (G) \Psi_0  ,\Psi_0 > .$$
Applying several times (\ref{Mizrahi-a}) with $X = e_j$,
we indeed obtain (\ref{prod-scal-base}). For every  $X$ in $H^2$, we have:
$$ \Big | (\partial _q - i \partial _{p})^{\alpha }  F(X) \Big |   \
 \Big |  ( \partial _q + i \partial _{p })^{\alpha } G(X) \Big |  \leq
  \Vert F \Vert _{Q }\
\Vert G \Vert _{Q } \ \prod _j ( Q (e_j, 0)^{1/2}  + Q(0, e_j)^{1/2}  ) ^{ 2\alpha _j}.  $$
 $$ \leq  2^{|\alpha|} \Vert F \Vert _{Q }\
\Vert G \Vert _{Q } \ \prod _j  ( Q (e_j, 0)  + Q(0, e_j)  ) ^{ \alpha _j}. $$
As a consequence,
$$ \sum _{|\alpha|= m} \left (\frac  {h^{|\alpha | }}{ 2^{|\alpha | }  \alpha !} \right )
\Big | (\partial _q - i \partial _{p})^{\alpha }  F (X)\Big |  \
\Big | ( \partial _q + i \partial _{p })^{\alpha } G(X) \Big | \leq
 \Vert F \Vert _{Q }\
\Vert G \Vert _{Q } \  \frac  {h^{m }}{   m  !}\left [ \sum _j  Q (e_j, 0)  + Q(0, e_j)  \right ] ^m $$
$$ \leq  \Vert F \Vert _{Q }\
\Vert G \Vert _{Q } \   \frac  {h^{m }}{   m  !} [ {\rm Tr } A_Q ]^m. $$
Therefore, the series (\ref{Miz-2}) converges absolutely for each  $X$.
According to these points,
$$ C_h^{wick} (F , G) (0) = \sum _{\alpha } <  Op_h^{wick} (G) \Psi_0 , u_{\alpha} >
 \ < u_{\alpha} ,  Op_h^{wick} (\overline F) \Psi_0 >  $$
 $$ = <  Op_h^{wick} (G) \Psi_0 ,  Op_h^{wick} (\overline F) \Psi_0 >
 =  <  Op_h^{wick} ( F) \circ Op_h^{wick} (G) \Psi_0 ,  \Psi_0 >
 = \sigma _h^{wick } \Big ( Op_h^{wick} ( F) \circ Op_h^{wick} (G) \Big ) (0). $$
 Consequently, equalities   (\ref{Miz-1}) and (\ref{Miz-2}) are proved for $X=0$. For any arbitrary other  $X$ one notes that, according to (\ref{alter-EC})
and (\ref{compos-Segal}),
 $$ \sigma _h^{wick } \Big ( Op_h^{wick} ( F) \circ Op_h^{wick} (G) \Big ) (X)=
 \sigma _h^{wick } \Big ( Op_h^{wick} ( F_X) \circ Op_h^{wick} (G_X) \Big ) (0)$$
where $F_X (U) = F(X+U)$.  Equalities (\ref{Miz-1}) and  (\ref{Miz-2}) then follow for 
any  $X$. The last inequality in the theorem  comes from the fact that:
 $$ \sum _{m+1}^{\infty}  \frac  {h^{k }}{   k  !} [ {\rm Tr } A_Q ]^k \leq
 h^{m+1}  [ {\rm Tr } A_Q ]^{m+1} e^{ h {\rm Tr } A_Q } .   $$

\bigskip

\subsection{Covariance by some Bogoliubov transformations.}\label{s-4-2}

\begin{theo}\label{t-Bogoliub} Let $U$ be a symplectic unitary linear map
in $H^2$. Then:

i) For each bounded operator $A$ in ${\cal F} _s (H_{\bf C})$, we have:
$$ \sigma _h^{wick} (  \Gamma (U)^{\star} A  \Gamma (U) ) (X) =
\sigma _h^{wick} (  A) ( UX). $$
ii) Let $\widetilde Q (X) = Q(UX)$. Then, for each $F\in S (H^2, Q)$, the
function $F \circ U$ belongs to $S (H^2, \widetilde Q)$, with the same norm.

iii) For each bounded operator $A$ in ${\cal L} (Q)$, the operator
$\Gamma (U)^{\star} A  \Gamma (U)$ is in   ${\cal L} (\widetilde Q)$,
with the same norm.

\end{theo}

Let us recall that, if a symplectic map $U$ in $H^2$ is not unitary, but verifying $U^{\star} U - I$ is trace class, then there is still a covariance for the Bogoliubov
transformation associated to $U$, but only for the Weyl calculus, not for the Wick
quantization. See \cite{A-L-N-3}.

{\it Proof.} If $U$ is both symplectic and unitary then it is also a ${\bf C}-$linear map
when $H^2$ is identified with $H_{\bf C}$. In other words, $U$ commutes with the map
$X \rightarrow \widehat X$.  Then, the coherent states
 $ \Psi_{X, h} $ (with $X\in H^2$) defined in (\ref{1.2}) satisfy:
$$ \Gamma (U)   \Psi_{X, h}    = \Psi_{UX, h}.
$$
Point i) follows easily. Point ii) follows from the chain rule. According to the above  point i) with  point ii)
of Theorem \ref{quantif-Wick} together with the facts that $U$ is a unitary operator commuting with the mapping $X \rightarrow \widehat X$, we have
$$ \Gamma (U) \Phi_S( \widehat V) \Gamma (U)^{\star} = \Phi_S(  \widehat {UV} ).
$$
As a consequence:
$$     [ \Phi_S( \widehat  V), \Gamma (U)^{\star} A  \Gamma (U) ] =
 \Gamma (U)^{\star} [ \Phi_S( \widehat {UV} ), A ]  \Gamma (U). $$
Iterating that process we obtain point iii).

\section{ Proof of Theorem  \ref{tPrincipal}. }\label{s-5}

\subsection{Proof of point i).}\label{s-5-3}

In the previous Sections \ref{s-2}, \ref{s-3} and \ref{s-4}, the space $H$ is an arbitrary infinite dimensional separable Hilbert space.
We now turn back to the photon Hilbert space of Section \ref{s-1}.  
We also consider the space   ${\cal H}_{sp}$ of Section \ref{s-1}. Let $Q$ be a quadratic form on $H^2$. The space $S(H^2, Q)$ of Definition \ref{d1.1} now refers as a space of functions taking values in ${\cal L}( {\cal H}_{sp})$ and the norm in the left hand side of  (\ref{1.1}) is the ${\cal L}( {\cal H}_{sp})$ norm. With a function $F$ in  $S(H^2, Q)$, we define an operator   $Op_h^{wick}(F)$
being bounded in  ${\cal H}_{ph} \otimes {\cal H}_{sp}$ where
${\cal H}_{ph} = {\cal F}_s (H_{\bf C})$.

Theorem \ref{t4-1} is then adapted.
Fixing a Hilbert basis  $(e_j)$ of $H$,
we set using the notations in Section \ref{s-4},  for all differentiable functions $F$ and $G$, with values in ${\cal L} ( {\cal H}_{sp})$,
$$ C^{1, wick } (F, G) = \frac {1} {2}
 \sum _{j}
 (\partial _{q_j} - i \partial _{p_j} )  F  \circ
( \partial _{q_j} + i \partial _{p_j } ) G , $$
where $\circ $ stands for the composition of operators in ${\cal L} ( {\cal H}_{sp})$.
Let $A$ be any operator written as   $A = \Phi_{S , h} (V) \otimes S $ where $V\in H^2$ and
$S\in {\cal L} ({\cal H} _{sp})$. Let $B= Op_h^{wick} (G)$ where
$G\in S(H^2 , Q)$. The Wick symbol of $A$ is the function 
$H^2 \ni X \rightarrow \varphi_V (X) = (V \cdot X) S$. The following equality comes from Theorem \ref{t4-1}:
\be\label{5-1} \sigma _h ^{wick} ( [A, B] ) (X) = [\varphi_V(X), G(X) ] +
h  C^{1, wick } (\varphi_V , G) (X) - h  C^{1, wick } (G, \varphi_V ) (X)   \ee
where the bracket in the first   term denotes the commutator of two operators
in   ${\cal L} ( {\cal H}_{sp})$.

We denote by ${\cal L} (Q)$ the space of bounded operators in ${\cal H}_{ph} \otimes {\cal H}_{sp}$
as in Definition \ref{d1.3} with $  e^{ \frac {i} {\sqrt {h}} \Phi_S ( t_1  \widehat  V_1+ ...
+ t_m  \widehat  V_m) }$ replaced by $  e^{ \frac {i} {\sqrt {h}} \Phi_S ( t_1  \widehat  V_1+ ...
+ t_m  \widehat  V_m) } \otimes I$.

For each $t\in \R$, we shall use the quadratic form $Q_t$ defined in
(\ref{def_Q_t}) and also use the quadratic form:
\be\label{tilde-Q-t} \widetilde Q_t (q, p) =  Q_t ( \chi_{-t} (q , p) ) =  |t|  \int_0^t
| d H^{free}_{int}(-s , q, p) |^2
ds.\ee

The following theorem is proved in \cite{A-L-N-3}.

 \begin{theo}\label{t8.1}
For any $t\in \R$, the family $   U_h^{red} (t)$ with any $h\in (0,1]$ defined by:
\be\label{Uh-red} U_h^{red} (t) =  e^{-i{t \over h}(H_{ph} \otimes I)}  e^{i{t \over h}   H(h) }\ee
 belongs to the class ${\cal L} (\widetilde Q_t)$  with the  quadratic form defined in
 (\ref{tilde-Q-t}).
Moreover,
\be\label{7.16} \Vert  U_h^{red} (t) \Vert _{{\cal L} (\widetilde Q_t, {\cal L} ({\cal H}_{sp}) )} = 1.\ee
 This function satisfies:
\be\label{Uh-red-equa}  \frac {d} {dt}  U_h^{red} (t) = i  H_{int}^{free} (-t)
U_h^{red} (t).
 \ee
 \end{theo}

  The point i) of Theorem \ref{tPrincipal}  will be a consequence of the following proposition.

\begin{prop}\label{p-classe} Let $A$ be a selfadjoint operator, bounded or unbounded in   ${\cal H} _{ph} \otimes {\cal H}_{sp}$. We suppose that $A$ is written under the form
 (\ref{forme-A})
where $F_A$ is an element of  $H^2$ and with
 $S_A$ lying in ${\cal L} ({\cal H}_{sp})$. Then:

 i) One has:
 $$  A^{free } ( t, h) = \Phi _S ( \chi_{-t}  (F_A)  ) \otimes I + I \otimes S_A. $$

 ii)  If $A(t, h)$ is defined in (\ref{A(t, h)}) then 
 the operator $A(t, h) - A^{free } ( t, h)$ is bounded in ${\cal L}(4Q_t)$  where $Q_t$ is the quadratic form defined in (\ref{def_Q_t}).

\end{prop}

{\it Proof.}  Point i) is a consequence of Theorem \ref{t-Bogoliub}. Let  $A$ be an operator as in (\ref{forme-A}).
 Set $A(t, h)$ the operator defined in (\ref{A(t, h)}).   Using the operator
$   U_h^{red} (t)$, $h\in (0,1]$, defined in (\ref{Uh-red}), one has:
$$ A(t, h) = e^{i{t \over h}(H_{ph} \otimes I)} U_h^{red} (t)  A
U_h^{red} (t)^{\star} e^{-i{t \over h}(H_{ph} \otimes I)}. $$
Since $U_h^{red} (t)$ is unitary:
$$ A(t, h) = e^{i{t \over h}(H_{ph} \otimes I)} \Big ( A +
[U_h^{red} (t) , A ] U_h^{red} (t)^{\star}  \Big )
e^{-i{t \over h}(H_{ph} \otimes I)}. $$
From Theorem \ref{t4-1}  and the expression (\ref{forme-A}) for $A$,
the commutator $[U_h^{red} (t) , A ]  $ belongs to ${\cal L}(\widetilde Q_t)$.
Consequently, by (\ref{comp}), the composition $[U_h^{red} (t) , A ]  U_h^{red} (t)^{\star} $ lies in ${\cal L}(4 \widetilde Q_t)$.
 According to Theorem \ref{t-Bogoliub}:
$$ e^{i{t \over h}(H_{ph} \otimes I)} [U_h^{red} (t) , A ] U_h^{red} (t)^{\star}
e^{-i{t \over h}(H_{ph} \otimes I)} \in {\cal L}(4 \widetilde Q_t \circ \chi_t).$$
Since $\widetilde Q_t \circ \chi_t = Q_t$, this proves Proposition \ref{p-classe}. \fpr

Set
$$ A^{[ red]} (t, h)  =  e^{-i{t \over h}(H_{ph} \otimes I)}
A(t, h) e^{i{t \over h}(H_{ph} \otimes I)}  =
 e^{-i{t \over h}(H_{ph} \otimes I)}  e^{i{t \over h}   H(h) } \
A\  e^{-i{t \over h}   H(h) }  e^{i{t \over h}(H_{ph} \otimes I)}
. $$
This operator satisfies:
\be\label{Ared(t)} \frac {d} {dt} A^{[red]}(t, h)  =
 i [H_{int} ^{free } (-t)   ,  A^{[red]}(t, h)],\qquad
  A^{[red]}(0, h)  = A. \ee
According to (\ref{5-1}):
\be\label{evol-2} \sigma _h^{wick} (  [H_{int} ^{free}(-t)  , A(t, h) ] ) (X) =
[ H_{int} ^{free}(-t, X , h) ,  A(t, X ,  h) ]
+ h C_1^{wick } ( H_{int} ^{free}(-t, \cdot  , h) ,  A(t, \cdot ,  h) ) (X)  \ee
$$ -  h C_1^{wick } (   A(t, \cdot ,  h)\  , \ H_{int} ^{free}(-t, \cdot  , h)   ) (X) $$
Consequently:
\be\label{equ-Ared}  {d\over dt}   A^{[  red]}  (t, X, h) =
 i [ H_{int} ^{free}(-t, X, h) ,  A^{[  red]}  (t, X, h) ] + \ee
 $$ + i C_1^{wick} \Big  (  H_{int} ^{free}(-t, \cdot ) \ ,
 \   A^{[ red]}  (t, \cdot , h )\Big ) (X) - i C_1^{wick} \Big  (  A^{[ red ]}  (t, \cdot , h ) \ , \
 H_{int} ^{free}(-t, \cdot ) \Big ) (X). $$

\subsection{Zeroth order term.}\label{s-A0}

The operator  $A^{[ red]} (t, h)$  satisfying  (\ref{Ared(t)}) is now approximated by a sum written under the following form:
$$ S^{[ M, red]}  (t,  h) = \sum _{j= 0}^M h^j A^{[ j, red]}(t)$$
where the operators  $A^{[ j, red]}(t)$ will be determined by their Wick symbols denoted $A^{[ j, red]}(t, X)$. 
Since the Wick symbol $ A^{[  red]}  (t, X, h)$ has to satisfy
 (\ref{equ-Ared}) then it makes sense to have in the case $j=0$:
\be\label{equ-tilde-A-0} {d\over dt}   A^{[ 0, red]}  (t, X) =
 i [ H_{int} ^{free}(-t, X) ,  A^{[ 0, red]}  (t, X) ]
 \ee
and $A^{[ 0, red]}  (0, X) = A(X)$, where $A(X)$ is the Wick symbol
of the observable $A$.  For $j\geq 1$, we must have:
\be\label{equ-tilde-A-j} {d\over dt}   A^{[ j, red]}  (t, X) =
 i [ H_{int} ^{free}(-t, X) ,  A^{[ j, red]}  (t, X) ] +
 \Phi ^{[ j-1, red ]}  (t,  X ) \ee
 where
 \be\label{equ-tilde-Phi-j} \Phi ^{[ j-1, red ]}  (t,  X ) =   i C_1^{wick} \Big  (  H_{int} ^{free}(-t, \cdot ) \ ,
 \   A^{[j-1, red]}  (t, \cdot , h )\Big ) (X) - i C_1^{wick} \Big  (  A^{[ j-1,red ]}  (t, \cdot , h ) \ , \
 H_{int} ^{free}(-t, \cdot ) \Big ) (X) \ee
and $A^{[ j, red]}  (0, X) = 0$.

In order to solve the equation (\ref{equ-tilde-A-0}), we use the following proposition,
which is  close to  Propositions 6.1 and 6.2 in \cite{A-L-N-2}.

\begin{prop}\label{pBloch1} There exists a unique function $G$
on $\R^2 \times H^2$  taking values in ${\cal L} ( {\cal H}_{sp})$ satisfying,
for all $X$ in $H^2$:
 \be\label{G(t)-1}{d\over dt} G(t, s,  X)
 = i G(t, s,  X)  H_{int }^{free} (t, X),\qquad
  H_{int }^{free} (t, X) =  H_{int } (\chi_t(X)),\ee
$$ G(s, s, X) = I.$$
Moreover,  $G(t, s, q, p)$ is unitary. If $0 \leq  s < t$ then the function
$X \rightarrow G(t, s, X)$ belongs to   $S(H^2 ,  Q_t)$ where $ Q_t$ is defined in
 (\ref{def_Q_t}) with a norm equal to 1.
This function also satisfies, for each $s$ and $t$ in $\R$, for each $X$ in $H^2$:
\be\label{groupe} G (s ,0,  X)^{\star} G(t, 0, X) = G( t-s , 0,  \chi_s (X)). \ee

\end{prop}

{\it Proof.} The existence property and the fact that $ G (t, s, X)$ is unitary can be considered as  standard facts. Next, we see that, for all $V$ in $H^2$:
 $$ {d \over dt } d  G(t, s,  X)(V) =
  i  G (t,  s,  X) d H_{int}^{free} (t, V) +
   i  d  G(t, s,  X)(V)  H^{free}_{int}(t, X). $$
 From Duhamel's principle and since $ d   G (s, s,   X)(V)
 = 0$ we then deduce:
 $$  d  G (t, s,  X)(V)  = i \int _s ^t
 G (\sigma , s,   X)
  d H^{free}_{int}(\sigma , V)
      G(t, \sigma ,  X)  d\sigma.   $$
Iterating, for all sequences $(V_1 , ... , V_m)$ in $H^2$:
 $$  d^m  G (t, s, \cdot)(V_1 , ... , V_m) =i^m
 \sum _{\varphi \in S_m} \int _{\Delta _m (t , s) }
   G ( \sigma _1 , s, \cdot )
  d  H^{free}_{int}(\sigma_1 , V_{\varphi (1)})
   G (\sigma _2, \sigma _{1} , \cdot )... $$
 $$
   ... d H^{free}_{int}(\sigma_{2} , V_{\varphi (2)})
   ...
   d H^{free}_{int}(\sigma_{m} , V_{\varphi (m)})
    G (t, \sigma _m,  \cdot ) d\sigma_1... d\sigma _m $$
where $S_m$ is the permutation group of $m$ elements and:
$$ \Delta _m (t , s) = \{ ( \sigma _1 , ... \sigma _m) , \ \ \ \
s < \sigma_1 < ... < \sigma _m < t \}.$$
Since $G(t, s, q, p)$ is unitary, one obtains:
$$ |d^m  G (t, s, \cdot)(V_1 , ... , V_m)| \leq
 \sum _{\varphi \in S_m} \int _{\Delta _m (t , s) }
 |d H^{free}_{int}(\sigma_1 , V_{\varphi (1)}) |...  | d H^{free}_{int}(\sigma_m , V_{\varphi (m)}) |
 d\sigma_1... d\sigma_m.$$
Consequently:
$$ | d^m  G (t, s, \cdot)(V_1 , ... , V_m) | \leq
 \prod _{j=1}^m \int _s^t
|d H^{free}_{int}(\sigma , V_j) | d\sigma
$$
where the above norm is the  ${\cal L} ( {\cal H}_{sp})$ norm. One notices that:
$$\int _0^t| d H^{free}_{int}(\sigma , V) |  d\sigma
\leq   Q_t (V )^{1/2}.$$
The first point then follows.
Concerning  equality (\ref{groupe}),
we remark that, as functions of the variable $t$, the two hand sides satisfy the same differential equation:
$$ {d\over dt} F(t, X) = i F(t, X)  H_{int }^{free} (t, X). $$
For the right hand  side of  (\ref{groupe}), we see that according to (\ref{wick-free}),
$ H_{int }^{free} (t-s,  \chi_s (X))=   H_{int }^{free} (t, X)$. Besides, the two hand sides of (\ref{groupe}) are equal for $t=s$. Therefore, they are equal everywhere.

If $G$ is the function of   Proposition \ref{pBloch1} and if $A$ is written as in (\ref{forme-A}) then
the function $A^{[ 0, red]} $ defined by:
\be\label{A-0-red} A^{[ 0, red]}  (t, X) = F_A (X)\otimes I +  G( -t, 0,  X)^{\star}  S_A  G(-t, 0, X) \ee
satisfies (\ref{equ-tilde-A-0}).
The first term, as a function of $X$,  is a continuous linear
function on $H^2$ and by  Proposition \ref{pBloch1}, the second term
belongs to $ S(H^2 , 4 \widetilde Q_t)$. The function defined by:
\be\label{A-0} A^{[ 0]}  (t, X) = F_A (\chi_t (X))\otimes I  + G( t, 0,  X) S_A G( t, 0,  X)^{\star}  \ee
satisfies  $ A^{[ 0, red]}  (t, X) =  A^{[ 0]}  (t, \chi_{-t} (X) )$ and this function
belongs to $ S(H^2 , 4  Q_t)$.

\subsection{$j$-th order term, $j\geq 1$.}\label{s-Aj}

Fix $j\geq 1$. Being given a function  $ A^{[ j-1, red]}  (t, X) $ in
$ S(H^2, 16^{j} \widetilde Q_t)$ we want to find a function
 $ A^{[ j, red]}  (t, X) $ satisfying (\ref{equ-tilde-A-j})
 with $\Phi ^{[ j-1, red ]}$ defined in (\ref{equ-tilde-Phi-j}).
  If $G$ is the function
of Proposition \ref{pBloch1} then the function defined, for each $X$ in $H^2$ by:
\be\label{A-j-red} A^{[ j, red]}  (t, X) = \int _0^t
 G(-t , 0, X)^{\star} G(-s , 0, X)
 \Phi ^{[ j-1, red]}  (s, X )
  G(-s , 0, X)^{\star} G(-t , 0, X)
  ds \ee
satisfies (\ref{equ-tilde-A-j})  and $A^{[ j, red]}  (0, X) = 0$. By (\ref{groupe}),
we also  have:
$$  A^{[ j, red]}  (t, X) = \int _0^t
G(t-s , 0, \chi_{-t}(X) )
 \Phi ^{[ j-1, red]}  (s, X )
 G(t-s , 0, \chi_{-t}(X) )^{\star} ds. $$
According to the induction hypothesis, if $0 \leq s \leq t$ then the  function $\Phi ^{[ j-1, red]}  (s, \cdot)$ belongs to
$S(H^2, 16^{j} \widetilde Q_s)$ also contained in $S(H^2, 16^{j} \widetilde Q_t)$. The function $X \rightarrow G(t-s , 0, \chi_{-t}(X) )$ is in $S(H^2,  \widetilde Q_t)$. By the result  (\ref{comp})  about composition,
the function $  A^{[ j, red]}  (t, \cdot) $ belongs to $S(H^2, 16^{j+1} \widetilde Q_t)$.
The function $  A^{[ j]}  (t, X)$ defined by:
\be\label{A-j} A^{[ j]}  (t, X) = \int _0^t
G(t-s , 0, X )
 \Phi ^{[ j-1, red]}  (s, \chi_t (X) )
 G(t-s , 0, X )^{\star} ds \ee
 where
\be\label{Phi-j} \Phi ^{[ j]}  (t,\cdot )
=i
C_1^{wick} ( H_{int} (\cdot) , A^{[ j-1]}  (t, \cdot))
-  i C_1^{wick} ( A^{[ j-1]}  (t, \cdot) , H_{int} (\cdot) )  \ee
satisfies $  A^{[ j, red]}  (t, X)=  A^{[ j]}  (t, \chi_{-t}(X))$.
This function is in $S(H^2, 16^{j+1} Q_t)$.

\subsection{End of the proof of Theorem  \ref{tPrincipal}.}\label{s-5-4}

The function  $ A^{[ 0]}(t, \cdot)$  of Theorem \ref{tPrincipal} is defined in (\ref{A-0}) and the functions  $ A^{[ j]}(t, \cdot)$ ($j\geq 1$) are constructed by iteration and given in  (\ref{A-j}) and  (\ref{Phi-j}).  We saw that $ A^{[ 0]}(t, \cdot)$ is the sum of a linear function on $H^2$ and of a function in
$S(H^2 , 4 Q_t)$. For $j\geq 1$,   $ A^{[ j]}(t, \cdot)$ belongs to $S(H^2 , 16^{j+1}  Q_t)$.
Therefore, the claims ii) and iii) of Theorem \ref{tPrincipal} are already proved for these functions.

   In order to prove (\ref{expansion}), we first derive a similar result for the
   functions  $ A^{[ j, red ]}(t, \cdot)$ defined in  (\ref{A-0-red})  for $j=0$,
   and in  (\ref{A-j-red}) and (\ref{equ-tilde-Phi-j}) for $j\geq 1$.
   For each integer $M$, let us set:
$$ S^{[ M, red]}  (t, X, h) = \sum _{j= 0}^M h^j A^{[ j, red]}(t, X).$$
By (\ref{A-0-red}), (\ref{A-j-red}) and (\ref{equ-tilde-Phi-j}), this function
satisfies:
$$  {d\over dt}   S^{[ M, red]}  (t, X, h) =
 i [ H_{int} ^{free}(-t, X) ,  S^{[ M, red]}  (t, X, h) ] + $$
 $$ + i h C_1^{wick} \Big  (  H_{int} ^{free}(-t, \cdot ) \ ,
 \   S^{[ M , red]}  (t, \cdot , h )\Big ) (X) - i h C_1^{wick} \Big  (  S^{[ M, red ]}  (t, \cdot , h ) \ , \
 H_{int} ^{free}(-t, \cdot ) \Big ) (X) $$
 $$ - h^{M+1} \Phi ^{[ M, red]} (t, X).$$
   We saw that $ A^{[ 0, red]}(t, \cdot)$ is the sum of a linear function on $H^2$ with a function in
$S(H^2 , 4 \widetilde Q_t)$.  For $j\geq 1$, $ A^{[ j, red]}(t, \cdot)$ belongs to $S(H^2 , 16^{j+1} \widetilde  Q_t)$. By theorem \ref{quantif-Wick}, there is
 an operator $ S^{[ M, red]}(t, h)$
 whose Wick symbol is $S^{[ M, red]}  (t, X, h)$. By Theorem \ref{t.Beals},
 this operator is the sum
of a Segal field with an operator in ${\cal L} (16^{M+1}  \widetilde Q_t)$.
Consequently, according to (\ref{5-1}):
 $$ {d\over dt}   S^{[ M, red]}  (t, h) =
  i [H_{int} ^{free } (-t)   ,  S^{[M, red]} (t, h)]
   - h^{M+1} Op_h^{wick}\Big ( \Phi ^{[ M, red]} (t, \cdot)\Big ). $$
Comparing with (\ref{Ared(t)}), we obtain:
$$ {d\over dt} \Big  (  A^{[red]}(t, h)   - S^{[ M, red]}  (t, h) \Big )  =
i [H_{int} ^{free } (-t)   , A^{[red]}(t, h)   -   S^{[M, red]}(t, h) ]  + $$
$$  + h^{M+1} Op_h^{wick} \Big ( \Phi ^{[ M, red]} (t, \cdot)\Big ).  $$
We also have $ A^{[ red ]}(0, h) - S^{[ M, red ]}(0, h)= 0$ and therefore, by Duhamel
principle:
$$  A^{[ red ]}(t, h) - S^{[ M, red]}(t, h) =  h^{M+1} R^{[ M, red ]}(t, h)$$
where
$$ R^{[ M, red ]}(t, h) = -  \int _0^t U_h^{red} (t) U_h^{red} (s)^{\star}
 Op_h^{wick} \Big ( \Phi ^{[ M, red]} (s, \cdot) \Big )  U_h^{red} (s) U_h^{red} (t)^{\star} ds, $$
 and $U_h^{red} (t)=e^{-i\frac{t}{h}H_{ph}}e^{i\frac{t}{h}H(h)}
$.
By Theorem \ref{t8.1} and since the quadratic form $\widetilde Q_t$ is an increasing
function of $t>0$, then the operators $U_h^{red} (t)$,
$U_h^{red} (s)$     ($0 \leq s \leq t$)  and their adjoints are in ${\cal L} (\widetilde Q_t)$.
We have seen that the function  $ \Phi ^{[ M, red]} (s, \cdot)$ belongs to  $S(H^2 , 16^{M+1}  \widetilde Q_s)$.
As a consequence, $Op_h^{wick} \Big ( \Phi ^{[ M, red]} (s, \cdot) \Big ) $ belongs to
${\cal L} ( 16^{M+1} \widetilde Q_s)$, thus in ${\cal L} ( 16^{M+1} \widetilde Q_t)$ if $0 \leq s \leq t$.
According to the property  (\ref{comp}) concerning the composition of operators in these classes, one observes that the operator $ R^{[M,  red ]}(t, h)$ is in
 ${\cal L} ( 16^{M+3} \widetilde Q_t)$.
According the covariance property, one therefore obtains that,
$$  A(t, h) - S^{[ M]}(t, h) =  h^{M+1} R^{[ M ]}(t, h)$$
where $ R^{[ M ]}(t, h)$ is in  ${\cal L} ( 16^{M+3}  Q_t)$.
Turning back to the Wick symbols,  point (\ref{expansion}) of Theorem 
\ref{tPrincipal} is then proved.

\section{Relation with physics equations (Maxwell and  Bloch).}\label{s-6}

We set, according to (\ref{A(t, h)}):
\be\label{9-1}B_j (x , t, h)=   e^{i{t \over h}   H(h) }( B_j (x )\otimes I)   e^{-i{t \over h} H(h)}  ,\qquad
E_j (x , t, h)=   e^{i{t \over h}   H(h) } (E_j (x )\otimes I)  e^{-i{t \over h} H(h)}, \ee
\be\label{9-2} S_j ^{[ \lambda]}  (t, h) =  e^{i{t \over h}   H(h) } ( I \otimes \sigma_j ^{[ \lambda]})   e^{-i{t \over h} H(h)}  . \ee

Since the initial observables  $B_j (x )\otimes I$, $E_j (x )\otimes I$ and  $I \otimes \sigma_j ^{[ \lambda]}$
are all under the form (\ref{forme-A}) then we can apply Theorem \ref{tPrincipal}.
The above Wick symbols have an asymptotic expansion in powers of $h$ with coefficients
defined by iteration in (\ref{A-j}). The aim of this section is
 to write this iteration under a form closer to the usual physics equations when the observables under consideration are the above observables. This is the content of Theorem \ref{tCalcul-termes} below.

To this end, we shall show that these operator valued functions satisfy the following equations which are similar to those in  Spohn \cite{SP}.
We use vectorial notations
 ${\bf B} (x, t, h) = ( B_1  (x , t, h) , B_2  (x , t, h) , B_3  (x , t, h) )$
and similarly for ${\bf E}  (x, t, h)$ and $ {\bf S}  ^{[\lambda ]} (t, h)$.
Given two operators triplets ${\bf A} = (A_1 , A_2 , A_3)$ and
${\bf B} = (B_1 , B_2 , B_3)$ we denote by ${\bf A} \times {\bf B}$ and
${\bf A} \times ^{sym} {\bf B}$ operators triplets defined by:
$$ \Big ( {\bf A} \times {\bf B} \Big ) _1 = A_2 B_3 - A_3 B_2,\qquad 
\Big ( {\bf A} \times ^{sym}  {\bf B} \Big ) _1 = \frac {1} {2}
( A_2 B_3 + B_3 A_2 - A_3 B_2 - B_2 A_3 ), $$
and the other components being similarly defined by circular permutations.
Thus,  one has $ {\bf A} \times ^{sym}  {\bf B}  = (1/2) (  {\bf A} \times {\bf B}  -
 {\bf B} \times {\bf A} )$.

 \begin{theo}\label{tSpohn} The operator valued functions defined in  (\ref{9-1}) and   (\ref{9-2}) satisfy, using the above vectorial notations:
\be\label{MaxOp-1} {\rm div} {\bf B}  (x, t, h) = {\rm div} {\bf E}  (x, t, h) = 0\ee
\be\label{Max-Op-2} {\partial \over \partial t}  {\bf B}    (x, t, h)  = - {\rm curl}
 {\bf E}   (x, t, h)\ee
\be\label{Max-Op-3} {\partial \over \partial t}  {\bf E}    (x, t, h)  =  {\rm curl}
 {\bf B}  (x, t, h)+ h \sum _{\lambda = 1} ^N {\bf S^{[\lambda ]}} (t, h) \times {\rm grad }
 \rho (x-x_{\lambda})\ee

where $\rho$ is defined in (\ref{rho}).
One also has:
\be\label{Bloch-Op} {d \over dt} {\bf S}  ^{[\lambda ]} (t, h) = 2 ( \beta + {\bf B}  (x_{\lambda} , t, h) )
\times ^{sym} {\bf S}  ^{[\lambda ]} (t, h).\ee
\end{theo}

{\it Proof.} Denoting, for instance, $B_j (x, q, p)$ as the Wick symbol  of the operator $B_j(x)$, we see that, for all $(q , p)$  in  $H^2$,
$$ \sum _{j= 1}^3 \frac {\partial B_j} {\partial x_j} (x , X) =
 \sum _{j= 1}^3 \frac {\partial E_j} {\partial x_j} (x , X) = 0.$$
One then deduces (\ref{MaxOp-1}), first for  $t=0$ and then for arbitrary $t$.
We also have, for  all $X$  in  $D(M)$
$$ \{ H_{ph}(\cdot)  ,  {\bf E}\}  (x , X) = {\rm curl }\  {\bf B} (x , X),\qquad
\{ H_{ph}(\cdot)  ,  {\bf B}\}  (x , X) = -  {\rm curl }\  {\bf E} (x , X). $$
Consequently, from Theorem \ref{t4.1}:
$$ [  H_{ph} ,  {\bf E} (x )] = (h/i) {\rm curl }\  {\bf B} (x) , \qquad
 [  H_{ph} ,  {\bf B} (x )] = - (h/i) {\rm curl }\  {\bf E} (x).$$
Using Definitions (\ref{7.3}), (\ref{7.4})  and (\ref{rho}),
we see that, for the Poisson brackets:
$$ \{ B_m (x , \cdot ), B_n (y , \cdot ) \} = 0 $$
$$ \{ E_m (x , \cdot ), B_n (y , \cdot ) \} =
{\rm grad \rho } (x - y) \cdot (e_m \times e_n ) $$
with $\rho$ defined in (\ref{rho}) and  where $(e_j)$ denotes the canonical basis of $\R^3$.
The Poisson bracket of two continuous linear forms  on  $H^2$ is independent of $(q, p)\in H^2$. It is here
a function depending only on  $x\in \R^3$. We then deduce, concerning the operators:
$$ [ E_m (x  ), B_n (y ) ] = (h/i) {\rm grad \rho } (x - y) \cdot (e_m \times e_n ). $$
Let us prove, for instance (\ref{Max-Op-3}). One has
$$ \frac {d} {dt} {\bf E}  (x , t, h) = (i/h)  e^{i(t/h) H(h)}
[ H(h) , {\bf E} (x)\otimes I ]  e^{- i(t/h) H(h)} $$
$$ = (i/h)  e^{i(t/h) H(h)}  \Big ( [H_{ph} , {\bf E} (x) ] \otimes I
\Big )  e^{-i(t/h) H(h)} +... $$
$$ ... + i \sum _{\mu = 1}^N \sum _{m=1} ^3   e^{i(t/h) H(h)}
 \Big (  [ B_m (x_{\mu}) , {\bf E} (x) ] \otimes  \sigma _m ^{[\mu ]}  \Big )
  e^{-i(t/h) H(h)}  $$
$$ =  e^{i(t/h) H(h)}  {\rm curl}\  {\bf B} (x) e^{-i(t/h) H(h)} +...$$
$$ ... - h  \sum _{\mu = 1}^N    e^{i(t/h) H(h)}
\Big  ( I \otimes {\rm grad}\ \rho (x - x_{\mu}) \times  \sigma  ^{[\mu ]}\Big )
 e^{-i(t/h) H(h)}.  $$
We then deduce (\ref{Max-Op-3}).
In order to prove   (\ref{Bloch-Op}), we see that, from (\ref{9-1}), (\ref{9-2}) and
(\ref{7.1}),  (\ref{7.8}):
$$ \frac {d} {dt}  S_j ^{[\lambda ]} (t, h) = (i/h)  e^{i(t/h) H(h)}
[ H(h),  I \otimes \sigma _j ^{[\lambda ]} ]  e^{i(t/h) H(h)} $$
$$ = i \sum _{\mu = 1}^N \sum _{m=1} ^3  e^{i(t/h) H(h)} (\beta_m+B_m ( x_{\mu}) )\otimes
 [ \sigma _m ^{[\mu ]} , \sigma _j ^{[\lambda ]}]  e^{-i(t/h) H(h)}.$$
One notices that $ [ \sigma _m ^{[\mu ]} , \sigma _j ^{[\lambda ]}] = 0$ if $\mu \not = \lambda$.
We then deduce that, for instance
$$ {d \over dt}  S_1 ^{[\lambda ]} (t, h) = 2  e^{i(t/h) H(h)}
\Big ( (\beta_2 +  B_2 ( x_{\lambda})) \otimes \sigma _3 ^{[\lambda ]} -
 (\beta_3 +B_3 ( x_{\lambda}))\otimes  \sigma _2 ^{[\lambda ]} \Big )
  e^{-i(t/h) H(h)}$$
$$ = 2 \Big ( (\beta_2 + B_2 ( x_{\lambda}, t, h) ) S_3 ^{[\lambda ]} (t, h) -
(\beta_3 +B_3 ( x_{\lambda}, t, h))  S_2 ^{[\lambda ]} (t, h) \Big ). $$
According to  (\ref{9-1}), (\ref{9-2}), the operators $B_j ( x, t, h)$
and  $S_k ^{[\lambda ]} (t, h)$ are commuting.
Thus one obtains   (\ref{Bloch-Op}). \fpr

Given two functions ${\bf F}$ and  ${\bf G}$,  $C^{\infty}$ on  $H^2$ and taking values in
$({\cal L} ({\cal H}_{sp}))$, we set, for each integer $j\geq 0$:
 $$ C^{j, wick} (F , G)  = \sum _{|\alpha| = j }
  \left (\frac  {1 } { 2^{|\alpha | }  \alpha !} \right )
(\partial _q - i \partial _{p})^{\alpha }  F  \circ
( \partial _q + i \partial _{p })^{\alpha } G
$$
where $\circ$ denotes the composition of two operators in ${\cal L} ({\cal H}_{sp})$.
Being given two  triplets ${\bf F}$ and  ${\bf G}$ of functions on $H^2$ taking values in $({\cal L} ({\cal H}_{sp}))$, we denote by
$ C ^{j, wick , \times}( {\bf F} , {\bf G})$ the function taking values in
$({\cal L} ({\cal H}_{sp}))^3$
defined by:
$$ \Big ( C ^{j, wick , \times}( {\bf F} , {\bf G}) \Big ) _1 =  C ^{j, wick } (F_2 , G_3) -
 C ^{j, wick } (F_3 , G_2),$$
the other component being similarly defined by circular permutations and set:
$$  C ^{j, wick , \times,sym}( {\bf F} , {\bf G})  =
\frac{1}{2}\Bigl(
C ^{j, wick , \times}( {\bf F} , {\bf G})-
C ^{j, wick , \times}( {\bf G} , {\bf F})
  \Bigr).$$

We denote by, for instance, $B_m (x, t, X, h)$ the  Wick symbol of the operator $B_m (x, t,  h)$
defined in  (\ref{9-1}).
Then, $B_m^{[j]}  (x, t, X)$ stands for the coefficient of $h^j$ in the asymptotic expansion of this symbol, given by Theorem \ref{tPrincipal}
(the operators considered here are all under the form (\ref{forme-A})).

\begin{theo}\label{tCalcul-termes}  With the above notations, one has, for  all  $j\geq 0$:
\be\label{MaxSym-1} {\rm div} {\bf B}^{[j]}  (x, t, X) = {\rm div} {\bf E}^{[j]}  (x, t, X) = 0\ee
\be\label{MaxSym-2} \frac  {\partial } { \partial t}  {\bf B}^{[j]}   (x, t, X)  = - {\rm curl}
 {\bf E}^{[j]}   (x, t, X).\ee
One has, for  $j=0$:
\be\label{MaxSym-3} {\partial \over \partial t}  {\bf E}^{[0]}    (x, t, X)  =  {\rm curl}
 {\bf B}^{[0]}  (x, t, X)) \ee
and  for  $j \geq 1$:
\be\label{MaxSym-4} {\partial \over \partial t}  {\bf E}^{[j]}    (x, t, X)  =  {\rm curl}
 {\bf B}^{[j]}  (x, t, X)+ \sum _{\lambda = 1} ^N {\bf S^{[\lambda, j-1 ]}} (t, X) \times {\rm grad }
 \rho (x-x_{\lambda}).\ee
One has for  $j=0$:
\be\label{MaxSym-5} {d \over dt} {\bf S}  ^{[\lambda, 0 ]} (t, X) = 2 ( \beta + {\bf B}^{[0]}  (x_{\lambda} , t, X) )
\times {\bf S}  ^{[\lambda , 0]} (t, X)\ee
and for  $j\geq 1$:
\be\label{MaxSym-6} {d \over dt} {\bf S}  ^{[\lambda , j ]} (t, X) = 2 ( \beta + {\bf B} ^{[0]}  (x_{\lambda} , t, X) )
\times {\bf S}  ^{[\lambda , j ]} (t, X) + ....\ee
$$ ... +  \sum _{ p+k+n =j , n <j}  C ^{p, wick , \times,sym}\Big (
{\bf B} ^{[k]}  (x_{\lambda} , t, \cdot ) \ ,\ {\bf S}  ^{[\lambda , n ]} (t, \cdot)\Big )(X).   $$
One has for  $t=0$ and  $j=0$:
$$ {\bf B}^{[0]}  (x, 0, X)  = {\bf B} (x, X) ,\qquad
{\bf E}^{[0]}  (x, 0, X)  = {\bf E} (x, X),\qquad
{\bf S}  ^{[\lambda , 0 ]} (0, X) = \sigma ^{[\lambda  ]}$$
and  for  $t=0$, $j\geq 1$:
$$ {\bf B}^{[j]}  (x, 0, X)  = 0,\qquad
{\bf E}^{[j]}  (x, 0, X)  = 0 ,\qquad
{\bf S}  ^{[\lambda , j ]} (0, X) = 0.$$

   \end{theo}

{\it Proof of equalities (\ref{MaxSym-1})-(\ref{MaxSym-4}).}
Let us prove, for example, equality  (\ref{MaxSym-4}).
Equality (\ref{Max-Op-3}) implies
\be\label{MaxSym-7}
   \frac  {\partial} {\partial t}  {\bf E}    (x, t, X, h)  =  {\rm curl}
 {\bf B}  (x, t, X, h)+ h \sum _{\lambda = 1} ^N {\bf S^{[\lambda ]}} (t, X, h) \times {\rm grad }
 \rho (x-x_{\lambda}). \ee
Let us show that the terms ${\bf E}^{[j]}    (x, t, X) $  are differentiable with respect to $t$ and that we have, in the sense of
Theorem \ref{tPrincipal}:
\be\label{9-6}  \frac  {\partial} {\partial t}  {\bf E}   (x, t, X, h) \sim
\sum _{r\geq 0} h^j \frac  {\partial} {\partial t}  {\bf E}^{[j]}    (x, t, X). \ee
To do this, we apply Theorem \ref{tPrincipal} with, for instance, the observable
$A = E_1(x)$ and also with the observable $B = (i/h) [H(h), A]$. The observable
$A$ is under the form (\ref{forme-A}). The operator  $B$ is also under this form since:
$$ B = (i/h) [H_{ph}, E_1 (x) ]\otimes I + i \sum _{\mu = 1}^N \sum _{m=1} ^3  [ B_m (x_{\mu}) ,
E_1 (x) ] \otimes \sigma_m^{[\mu ]} $$
$$ = \left ( \frac {\partial B_3(x) } {\partial x_2} - \frac {\partial B_2(x) } {\partial x_3}
\right )\otimes I
- h \sum _{\mu = 1}^N \sum _{m=1} ^3 ({\rm grad} \rho ( x - x_{\mu}) \cdot ( e_1 \times e_m) )
(I \otimes \sigma_m^{[\mu ]}). $$
According to  Theorem \ref{tPrincipal}, the Wick symbol  of the operator
$ (i/h) e^{i(t/h) H(h)}[H(h), E_1 (x)] e^{-i(t/h) H(h)} $ has an asymptotic expansion
described in this theorem. The coefficient of $h^j$ in this development is the derivative of $E_1^{ [j]} (x , t, q, p)$ and  we indeed have (\ref{9-6}). We similarly prove that, in the sense of Theorem \ref{tPrincipal}:
$$  {\rm curl} \   {\bf B}   (x, t, X, h) \sim
\sum _{j\geq 0} h^j    {\bf B}^{[j]}    (x, t, X). $$
 Consequently  the two hand sides in (\ref{MaxSym-7}) have asymptotic expansions in powers of $h$. Identifying the coefficients of $h^j$ in the two hand sides,
 we then deduce  (\ref{MaxSym-4}).

 {\it Proof of equalities (\ref{MaxSym-5}) and   (\ref{MaxSym-6}).}
 Equality (\ref{Bloch-Op}) implies:
 $$ \frac  {d } { dt} {\bf S}  ^{[\lambda ]} (t, X, h) =  \sigma _h ^{wick} \Big ( ( \beta + {\bf B}  (x_{\lambda} , t, X, h) ) \times {\bf S}  ^{[\lambda ]} (t,X, h) \Big )
 + ... $$
 $$ ... -  \sigma _h ^{wick} \Big (  {\bf S}  ^{[\lambda ]} (t, h)  \times
 ( \beta + {\bf B}  (x_{\lambda} , t,  h) ) \Big )(X)
 .$$
From Theorem  \ref{tPrincipal}, we have, for  all integers  $M$:
 $$  \sigma _h ^{wick}  ( {\bf B}  (x_{\lambda} , t) ) (X)  = \sum _{k= 0}^M
 h^k {\bf B}^{[k]}  (x_{\lambda} , t, q, p) + h^{M+1} R_M (t, q, p, h ) $$
 $$  \sigma _h ^{wick}  (  {\bf S}  ^{[\lambda ]} (t,h) )=  \sum _{n= 0}^M
 h^n {\bf S}^{[\lambda , n]}  (x_{\lambda} , t, X) + h^{M+1} S_M (t, X, h ) $$
 where    $R_M (t, \cdot, h)$ and  $S_M (t, \cdot, h)$ belong to $S(H^2, 16^{M+5} Q_t)$,
 with a norm  bounded
uniformly in $t$ and  $h$ when $t$ belongs to a compact set of $\R$ and  when $h$
belongs to $(0, 1]$.
In view of Theorem   \ref{t-Mizrahi},  we then deduce:
 \be\label{6-XX} \sigma _h ^{wick} \Big ( ( \beta + {\bf B}  (x_{\lambda} , t, h) ) \times {\bf S}
  ^{[\lambda ]} (t, h) \Big )  =  \sum _{j= 0}^M h^j  ( \beta + {\bf B}^{[0]}  (x_{\lambda} , t, \cdot) )
\times {\bf S}  ^{[\lambda , j]} (t, \cdot ) + ...\ee 
  $$...  +  \sum _{p+ k + n \leq M, p+ k>0} h^{p+k+n}  C ^{p, wick , \times,sym}\Big (
{\bf B} ^{[k]}  (x_{\lambda} , t, \cdot ) \ ,\ {\bf S}  ^{[\lambda , n ]} (t, \cdot)\Big )
+ h^{M+1} T_M (t,\cdot , h )   $$
where  $T_M (t,\cdot , h )$ belongs to some space  $S(H^2, K Q_t)$ with some $K>0$ and
with a norm bounded uniformly in $t$ and  $h$ when $t$ belongs to a compact set of  $\R$ and when $h$
is in $(0, 1]$. As above, we show that ${\bf S}  ^{[\lambda , j ]} (t, X )$ is differentiable with respect to
$t$ and that we have, in the sense of Theorem \ref{tPrincipal},
\be\label{MaxSym-8}  \frac  {d } { dt} {\bf S}  ^{[\lambda ]} (t, X, h)  \sim \sum _{j\geq 0} h^j
  \frac  {d } { dt} {\bf S}  ^{[\lambda , j ]} (t, X). \ee
When identifying the coefficients of $h^j$ in (\ref{6-XX})  and  in
 (\ref{MaxSym-8}), one then obtains    (\ref{MaxSym-5}) and   (\ref{MaxSym-6}). \fpr

For example, the first term satisfies the  Bloch equations:
\be\label{9-3} \frac  {d } { dt} {\bf S}  ^{[\lambda , 0 ]} (t, X) =
2 ( \beta + {\bf B} ^{[0]}  (x_{\lambda} , t, X) )
\times {\bf S}  ^{[\lambda , 0 ]} (t, X). \ee
The second term satisfies, from (\ref{MaxSym-6}) and Theorem \ref{t-Mizrahi}: 
\be\label{9-4} \frac  {d } { dt} {\bf S}  ^{[\lambda , 1 ]} (t, X) =
2 ( \beta + {\bf B} ^{[0]}  (x_{\lambda} , t, X) )
\times {\bf S}  ^{[\lambda , 1 ]} (t, X) + \cdots \ee
$$ \cdots +  2
{\bf B} ^{[1]}  (x_{\lambda} , t, X) \times  {\bf S}  ^{[\lambda , 0 ]} (t, X)
 + {\bf K}  ^{[\lambda , 1 ]} (t, X)  $$
with, for instance:
\be\label{9-5}  K_1  ^{[\lambda , 1 ]} (t,\cdot )  = d B_2 ^{[0]} (x_{\lambda} , t, \cdot )
\cdot  d S_3 ^{[\lambda , 0]} ( t, \cdot )
- d B_3 ^{[0]} (x_{\lambda} , t, \cdot )
\cdot  d S_2 ^{[\lambda , 0]} ( t, \cdot ).\ee
We have denoted by $(dF) (q , p) \cdot (dG) (q , p)$  the scalar product of the differentials of  two functions on $H^2$.
The second term in the right hand side of (\ref{9-4}) only reflects the influence, according to the classical Bloch equations, of the radiated field, according to Maxwell equations, by all the spins
between times  $0$ and  $t$. Only  the term $ K_1  ^{[\lambda , 1 ]} (t,\cdot ) $ is genuinely a quantum correction.

\section{Photon emission semiclassical study.}\label{s-7}

This section is concerned with the proof of Theorem \ref{t-evol-N}. To this end, the operators    $E_{j}^{pol} (x )$ involved in this result are now precisely defined.

We need to introduce, at each point $x$ in $\R^3$, six  operators
 $B_j^{pol} (x)$ and $E_j^{pol} (x)$
($1\leq j \leq 3$)  having no any counterpart in classical physics.
We set ${\cal F} (q , p) = (-p, q)$.  We denote by $E_+$
the subspace of all   $X\in H^2$ satisfying $JX = {\cal F} X$, where
$J$ is defined in (\ref{7.5}),   and by $E_-$
the subspace of all
$X\in H^2$ verifying $JX = - {\cal F} X$. These two subspaces correspond
to the circular right and left polarization notions in physics.
Then,  $\Pi_{\pm}: H^2 \rightarrow  E_{\pm}$ stands for the
corresponding orthogonal projections. One has,
$$ \Pi_+X = {1\over 2} \Big ( X - J {\cal F} X \Big ),\quad
\Pi_-X = {1\over 2} \Big ( X +  J {\cal F} X \Big ),$$
and thus, $\Pi_+ - \Pi_-= - J {\cal F}$. Set,
 \be\label{10-2}  B_{j}^{pol} (x , q, p) = -  B_{j} (x ,J {\cal F} (q , p)),
\quad E_{j}^{pol} (x , q, p) = -  E_{j} (x ,J {\cal F} (q , p)).  \ee
We denote by
 $  B_{j}^{pol}(x)$  and $  E_{j}^{pol}(x)$ the operators whose Wick symbols are
 $  B_{j}^{pol} (x, q, p)$ and $  E_{j}^{pol} (x, q, p)$.  We set:
\be\label{10-3} B_j ^{pol}(x , t, h)=   e^{i{t \over h}   H(h) }
( B_j^{pol} (x) \otimes I)
  e^{-i{t \over h} H(h)},\qquad
E_j ^{pol}(x , t, h)=   e^{i{t \over h}   H(h) }
( E_j^{pol} (x) \otimes I)
  e^{-i{t \over h} H(h)}
  .\ee

\begin{lemm}\label{t5.1}
The operator $N'( t, h)$ defined in  (\ref{N'(t)}) satisfies:
\be\label{equation-1-N}  N'( t, h)=  \sum _{\lambda =1}^N  \sum _{m=1}^3
 e^{i{t \over h}   H(h) } (  E_m^{pol} (x_{\lambda})\otimes
 \sigma_m^{[\lambda]} ) e^{-i{t \over h}   H(h) }. \ee
This operator is bounded from $D(H(h))$ to ${\cal F}_s(H_{\bf C})\otimes {\cal H}_{ph}$. We have also:
\be\label{equation-2-N}  N'( t, h) =   -\sum _{\lambda = 1}^N \sum _{m=1}^3   E_m^{pol} (x_{\lambda }, t, h)
 \circ  S_m  ^{[\lambda ]} (t , h). \ee
\end{lemm}

{\it Proof.} We begin with (\ref{N'(t)}).
Clearly, one has $ [ H_{ph} \otimes I , N \otimes I ] = 0$ since $H_{ph} = h {\rm d} \Gamma (M_{\omega})$,
$N = {\rm d} \Gamma (I)$ and since $M_{\omega}$ commutes with $I$. Therefore:
$$ [H_{int} , N \otimes I ] =  \sum _{\lambda =1}^N  \sum _{m=1}^3
[ B_m (x_{\lambda}), N ] \otimes
\sigma_m^{[\lambda]}.$$
By \cite{DG} (Lemma 2.5 iii), third identity), we have:
$$ i  [ B_m (x), N ] =  - E_m^{pol} (x),$$
for all $x\in \R^3$. This point comes from standard commutations properties (see {\it e.g.}, Lemma 2.5
(ii) in \cite{DG}). Therefore, (\ref{equation-1-N}) follows from (\ref{N'(t)}) and these previous computations.
We saw that the Segal fields $ E_m^{pol} (x_{\lambda}) $ are bounded from $D(H_{ph}) $ to  ${\cal F}_s(H_{\bf C})$,
because the corresponding elements of $H^2$ are in $D(M_{\omega} ^{-1/2})$. Moreover the domain
of $H(h)$ is $D( H_{ph}) \otimes {\cal H}_{sp}$. The equality (\ref{equation-2-N}) is another
formulation of (\ref{equation-1-N}).    \fpr

{\it End of the proof of Theorem  \ref{t-evol-N}.}
The operators $  B_{j}^{pol}(x)$  and  $  E_{j}^{pol}(x)$ defined in (\ref{10-2})
are under the form (\ref{forme-A}). Theorem \ref{tPrincipal} shows that each
operator $ B_j ^{pol}(x , t, h) $ and  $ E_j ^{pol}(x , t, h) $ defined in (\ref{10-3}) is the sum  of
a Segal field with an operator belonging to ${\cal L}(4Q_t)$, in the sense of Section
\ref{s-5}.   Theorem \ref{tPrincipal} shows that the symbols have asymptotic expansions that can be written, in order to simplify notations, as:
$$ \sigma_h^{wick} ( {\bf E}^{pol} (x, t, h))(X) = {\bf E}^{pol} (x, t, X,  h)
 \sim \sum _{j \geq 0}
h^j {\bf E}^{pol, j} (x, t, X).$$
Then Theorem \ref{t-evol-N} follows from (\ref{equation-2-N}), from Theorems \ref{t4-1}  and
\ref{t-Mizrahi}, and from the above asymptotic expansions. We also have:
$$ N^{[j]}  (t, \cdot) = \sum _{\lambda = 1}^N \sum _{q=1}^3 \sum _{k+m+n = j}
C^{k, wick}  \Big (  E_q^{pol, m} (x_{\lambda}, t, \cdot ) \ , \
S_q^{[\lambda , n ]}  ( t, q, p ) \Big ). $$

\fpr

laurent.amour@univ-reims.fr\newline
LMR EA 4535 and FR CNRS 3399, Universit\'e de Reims Champagne-Ardenne,
 Moulin de la Housse, BP 1039,
 51687 REIMS Cedex 2, France.

lisette.jager@univ-reims.fr\newline
LMR EA 4535 and FR CNRS 3399, Universit\'e de Reims Champagne-Ardenne,
 Moulin de la Housse, BP 1039,
 51687 REIMS Cedex 2, France.

jean.nourrigat@univ-reims.fr\newline
LMR EA 4535 and FR CNRS 3399, Universit\'e de Reims Champagne-Ardenne,
 Moulin de la Housse, BP 1039,
 51687 REIMS Cedex 2, France.

\end{document}